\def\bibart#1#2#3#4#5#6#7
{\bibitem{#1}
#2, #4, #5 #6 (#3), #7}

\def\bibcoll#1#2#3#4#5#6#7#8
{\bibitem{#1}
#2, #4, in: #5, #6, #7, #3, pp. #8}

\def\bibbook#1#2#3#4#5#6
{\bibitem{#1} #2, #4, #5, #6, #3}

\def\bibdiss#1#2#3#4#5#6
{\bibitem{#1}
#2 (#3), #4, #5, #6}

\def\mat#1{\mathbb{#1}}

\def\qed{\hfill\ \rule{2mm}{2mm} }

\def\qex{\hfill\ \vbox{\hrule\hbox{\vrule\kern4pt\vbox{\kern4pt{}
\kern4pt}\kern4pt\vrule}\hrule}}

\def\hochdamit{\! \Rsh }

\def\bothways{\matrix{\vartriangleright \vspace*{-.1in} \cr \vartriangleleft \cr } }

\def\directbothways{\matrix{\blacktriangleright \vspace*{-.1in} \cr \blacktriangleleft \cr } }

\def\transcont{\vartriangleright \hspace*{-.08in} \vartriangleright }

\def\setcontrol{\searrow }

\newcounter{casectr}

\newcounter{claimctr}

\documentclass[12pt]{article}

\usepackage{amssymb}
\newtheorem{guess}{Guess}[section]

\newtheorem{define}[guess]{Definition}
\newtheorem{prop}[guess]{Proposition}
\newtheorem{theorem}[guess]{Theorem}

\newtheorem{cor}[guess]{Corollary}

\newtheorem{lem}[guess]{Lemma}

\newtheorem{nota}[guess]{Notation}

\newtheorem{remark}[guess]{Remark}

\newtheorem{exam}[guess]{Example}

\newtheorem{oques}[guess]{Open Question}

\usepackage[pdftitle={Automorphism Conjecture}, pdfauthor={Bernd Schroeder},
colorlinks=true, pdfstartview=FitV, linkcolor=blue, citecolor=blue, urlcolor=blue]{hyperref}

\usepackage{color}

\begin{document}

\bibliographystyle{plain}

\title{
Towards the Automorphism Conjecture I:
Combinatorial Control
}

\author{
\small Bernd S. W. Schr\"oder \\
\small School of Mathematics and Natural Sciences\\
\small The University of Southern Mississippi\\
\small
118 College Avenue, \#5043\\
\small Hattiesburg, MS 39406\\
}

\date{\small \today}

\maketitle

\begin{abstract}

This paper
exploits
adjacencies between
the orbits of an ordered set $P$
and a consequence of the classification of finite simple groups
to, in many cases,
bound
the number of automorphisms
by
$2^{0.99|P|} $.
Results clearly identify the structures
which currently prevent the proof of
such an exponential bound, or
which
indeed
inflate the number of automorphisms
beyond such a bound.
This is a first step towards
a possible resolution of the
Automorphism
Conjecture for ordered sets.

\end{abstract}


\section{Introduction}

An {\bf ordered set}
consists of an underlying set $P$ equipped with a
reflexive, antisymmetric and transitive relation $\leq $, the order
relation.
An {\bf order-preserving self-map},
or, an {\bf endomorphism},
of an ordered set $P$ is a
self map
$f:P\to P$ such that $p\leq q$ implies
$f(p)\leq f(q)$.
Consistent with standard terminology,
endomorphisms with an inverse that is
an endomorphism, too, are called {\bf automorphisms}.
The set of endomorphisms is denoted ${\rm End} (P)$
and the set of automorphisms is denoted ${\rm Aut} (P)$.
Rival and Rutkowski's
Automorphism Problem (see \cite{RiRut}, Problem 3)
asks the following.

\begin{oques}

{\bf (Automorphism Problem.)}
Is it true that
$$\lim _{n\to \infty } \max _{|P|=n}
{ |{\rm Aut} (P)|\over |{\rm End} (P)| } =0?$$

\end{oques}

The {\bf Automorphism Conjecture} states that
the
Automorphism Problem
has an affirmative answer.
In light of the facts that, for
almost every ordered set, the identity is the only
automorphism
(see \cite{Proe}, Corollary 2.3a), that
every ordered set
with $n$ elements
has at least $2^{{2\over3} n} $ endomorphisms
(see \cite{DRSW}, Theorem 1),
and that this lower bound has been improved to
at least $2^{n} $ endomorphisms
in the unpublished manuscript
\cite{Duffusprivate},
this conjecture
is quite natural.
However, the Automorphism Conjecture has been remarkably resilient
against attempts to prove it in general.

With lower bounds on the number of endomorphisms available,
it is natural to also consider upper bounds on the number of
automorphisms.
Section \ref{ACnoprobsec}
shows that,
for the Automorphism Conjecture,
we need not be concerned with
order-autonomous antichains.
Section \ref{primnestsec}
presents
key results from the
theory of permutation groups
and
Section \ref{controlsetsec}
shows
how to use them to
recursively generate an upper bound for the number of
automorphisms.
Recall that a {\bf block} for a permutation group is a set $B$ that,
under each permutation, is either invariant or it is mapped to
a set disjoint from $B$.
Section \ref{combcontsec}
presents the key notion of
combinatorial control between
blocks
of groups of automorphisms.
Section \ref{depth0sec}
shows that
block sets on which the automorphism group induces an
alternating or a symmetric group are
indeed the only obstacles to
obtaining
a bound of the form
$2^{0.99|P|} $.
Section \ref{factpairinteract}
presents a key notion in which such block sets can be
interrelated.
Section \ref{maxlocksec} shows how to
compensate for factorials when the
connected components with respect
to this interrelation become too large.
Section \ref{mainres} provides the main result
with explicit references to all necessary definitions
and Section \ref{takestock}
provides a quick application
to ordered sets of width at most $9$.

\section{Exempting Order-Autonomous Antichains
from Further Analysis}
\label{ACnoprobsec}

Let $P$ be an ordered set.
For sets $B,T\subseteq P$, we will write
$B<T$ iff every $b\in B$ is strictly below every $t\in T$.
For singleton sets, we will omit the set braces.
Recall that a nonempty subset $A\subseteq P$ is called
{\bf order-autonomous}
iff, for all $z\in P\setminus A$, we have that
existence of an $a\in A$ with $z<a$ implies $z<A$, and,
existence of an $a\in A$ with $z>a$ implies $z>A$.
An
order-autonomous subset $A\subseteq P$ will be called {\bf nontrivial}
iff $|A|\not\in \{ 1,|P|\} $.
Further recall that an {\bf antichain} is a subset in which no two distinct elements
are comparable.

Every permutation of a
nontrivial order-autonomous antichain can be
extended to an automorphism of the surrounding ordered set
by mapping all other points to themselves.
Therefore,
nontrivial order-autonomous antichains
inflate the
number of automorphisms and they
tend to clutter
any estimates of their number.
Proposition \ref{onelexiter} below shows that,
if the Automorphism Conjecture holds
in the absence of nontrivial order-autonomous antichains,
then the Automorphism Conjecture also holds in the presence of
nontrivial order-autonomous antichains.
The underlying idea is a
technique used, for  example, in \cite{LRZ,LiuWan}:
When there are ``many" automorphisms,
we show that
there are
``enough corresponding endomorphisms" to guarantee the
automorphism to endomorphism
ratio's convergence to zero.

\begin{prop}
\label{onelexiter}

Let ${\cal T}$ be a class of finite
ordered sets such that no set in
${\cal T}$ has a nontrivial order-autonomous antichain
and for which the Automorphism
Conjecture holds.
Let
${\cal A} $ be the
class of
ordered sets that can be obtained from
ordered sets $T\in {\cal T}$ by replacing every
element of $T$ with a nonempty, but possibly trivial,
order-autonomous antichain.
Then
the Automorphism
Conjecture holds for
${\cal A} $.

\end{prop}

{\bf Proof.}
Let
$\varepsilon >0$.
Fix
$K\in {\mat N}$ such that
${
K!
\over
K^K
}
<\varepsilon $
and
$2^{-K} <\varepsilon $.
For every $n\in {\mat N}$, let
$h(n):=
\sup _{T\in {\cal T} , |T|\geq n}
{ |{\rm Aut} (T)|\over |{\rm End} (T)| } .
$
Choose $N\in {\mat N}$ such that, for all
$n\geq N$, we have
$h\left( n-
K^2 \right) \leq {\varepsilon \over \left( K^2 \right) !} $.
Now let
$P\in {\cal A}$
have $n\geq N$ elements.

Let $M\subseteq {\rm Aut} (P)$ be the
subgroup of automorphisms that
map every maximal
(with respect to inclusion)
order-autonomous antichain
to itself.
Let $\Phi \in M$.
Because every $\Psi \in {\rm Aut} (P)$ maps
maximal order-autonomous antichains
to maximal
order-autonomous antichains,
every
$\Psi ^{-1} \Phi \Psi $
maps every
maximal
order-autonomous antichain
to itself.
Therefore $M$ is a normal subgroup of ${\rm Aut} (P)$
and we have
$|{\rm Aut} (P)|=|{\rm Aut} (P)/M|\cdot |M|$.

Let $m_1 , \ldots , m_\ell $ be the individual sizes of the
nontrivial
maximal order-autonomous antichains
in $P$.
Then $|M|=
\prod _{j=1} ^\ell m_j !$.

{\em Case 1:
$\sum _{j=1} ^\ell m_j >K^2 $.}
Every order-autonomous antichain
with $m_j $ elements has a total of
$m_j ^{m_j } $ self-maps.
Hence, in every $[\Psi ]\in {\rm Aut} (P)/M$,
there are exactly
$\prod _{j=1} ^k m_j !$
automorphisms of $P$ and,
for every $[\Psi ]\in {\rm Aut} (P)/M$,
there are at least
$\prod _{j=1} ^k m_j ^{m_j } $
corresponding endomorphisms.
Hence
${|{\rm Aut}
(P)|
\over
|{\rm End}
(P)|}
\leq
{|{\rm Aut} (P)/M|
\prod _{j=1} ^k m_j !
\over
|{\rm Aut} (P)/M|
\prod _{j=1} ^k m_j ^{m_j }
}
=
\prod _{j=1} ^k {m_j !\over m_j ^{m_j } }
$.
Now,
if there is an
$m_j \geq K$,
we obtain
$
{
|{\rm Aut} (P )|
\over
|{\rm End} (P)|
}
\leq
{
m_j !
\over
m_j ^{m_j }
}
\leq
{
K!
\over
K^K
}
<\varepsilon
.$
Otherwise, for all $j\in \{ 1, \ldots , \ell \} $, we
have
$m_j <K$.
In this case,
$\ell \geq K$
and
hence
${
|{\rm Aut} (P)|
\over
|{\rm End} (P)|
}
\leq
\prod _{j=1} ^\ell
{m_j !
\over
m_j ^{m_j }
}
\leq
\left(
{
2
\over
4
}
\right) ^\ell
\leq
\left(
{
1
\over
2
}
\right) ^K
<\varepsilon
.
$

{\em Case 2:
$\sum _{j=1} ^\ell m_j \leq K^2 $.}
Let $T\subseteq P$ be obtained by,
from every nontrivial maximal
order-autonomous antichain,
deleting all but one
element.
For every $t\in T$, let $A_t $ be the unique
maximal order-autonomous antichain that contains $t$.
Then,
for every
equivalence class $[\Psi ]\in {\rm Aut} (P)/M$,
and all
$\Phi \in [\Psi ]$, we have that
$\Phi [A_t ]\cap T=\Psi [A_t ]\cap T$.
Hence if, for every $[\Psi ]\in {\rm Aut} (P)/M$ and
every $t\in T$, we set
$\Psi _T [t]:=\Psi [A_t ]\cap T$, we obtain a well-defined function, and
it is easy to check that $\Psi _T \in {\rm Aut} (T)$.
Via the map $[\Psi ]\mapsto \Psi _T $,
we see that
${\rm Aut} (T)$
contains a subgroup that is isomorphic to
the factor group
${\rm Aut} (P)/M$.
Moreover, because every endomorphism $f$ of $T$
can be extended to an endomorphism of $P$
by mapping every $A_t $ to $f(t)$, we have
$|{\rm End} (P )|\geq |{\rm End} (T )|$.
Therefore
$
{
|{\rm Aut} (P )|
\over
|{\rm End} (P )|
}
\leq
{
|{\rm Aut} (T)|\cdot \left( \prod _{j=1} ^\ell m_j !\right)
\over
|{\rm End} (T )|
}
\leq
h\left(
|T|
\right)
\left( \sum _{j=1} ^\ell m_j \right) !
\leq
h\left(
n
-
K^2
\right)
\left( K^2 \right) !
<\varepsilon
.$
\qed

\section{Primitive Nestings}
\label{primnestsec}

Because the automorphism group is a special type of permutation
group, we need
key concepts and
deep results from the theory of permutation groups.
In the case of Corollary \ref{marotiexpbound} below,
a
simple computation leads to a simple improvement
of a result by A. Mar\'oti.
As indicated in the Mathematical Review MR1943938
of \cite{Maroti},
Theorem \ref{24special}
below rests on the shoulders of the giant which is the
Classification of Finite Simple Groups.
Consequently, the subsequent
results on the Automorphism Problem
are rather deep by default.
Thanks to the accessible writing in \cite{BuekLee}
and \cite{Maroti}, this depth is
within reach of the author's humble combinatorial means.

We should note that
the results from \cite{Maroti}
are also used to compute
upper bounds for the order of
automorphism groups of graphs, see,
for example, \cite{BabaiSRGAut}.
However, our goal here
is an exponential bound of $2^{0.99n} $
plus a clear understanding
of the structure of the sets for which this bound does not hold,
see Theorem \ref{upperboundonlyoutsizedleft} below.
Meanwhile, for graphs,
the goal
is a subexponential bound, which necessitates a focus on
very special classes of graphs.
Consequently, the approach here is complementary to
corresponding work in graph theory, with the intersection
likely being the use of \cite{Maroti}.

\begin{define}

Let $G$ be a permutation group on the set $X$.
The {\bf degree} of $G$
is the size $|X|$ of the underlying set $X$.
For every $x\in X$, we define
$G\cdot x:=\{ \sigma (x):\sigma \in G\} $ and call it the
{\bf ($G$-)orbit} of $x$.
We call $G$ {\bf transitive}
iff, for one (and hence for all) $x\in X$, we have
$G\cdot x=X$.

\end{define}

Without further knowledge about interactions between orbits,
it is natural to focus on transitive permutation groups.
However, the requisite definitions and many elementary results
do not require transitivity.
Hence, in the following, transitivity is only needed when
explicitly stated or part of the definition.
Starting in Section \ref{combcontsec}, the interactions between
blocks in distinct orbits will take center stage.

\begin{define}
\label{blockdef}

Let $G$ be a permutation group on the set $X$.
A subset $B\subseteq X$ is called a {\bf ($G$-)block}
iff, for all $\sigma \in G$, we have
%that
$\sigma [B]=B$ or $\sigma [B]\cap B=\emptyset $.
A block is called {\bf nontrivial} iff it is not
a singleton and not equal to
$X$.
A transitive permutation group $G$ is called {\bf primitive}
iff it has no nontrivial $G$-blocks.

\end{define}

\begin{theorem}
\label{24special}

(See Corollary 1.4 in \cite{Maroti}.)
Let $G$ be a
primitive subgroup of the symmetric group $S_n $.
that does not contain
the alternating group
$A_n $. If $|G|>2^{n-1} $, then $G$ has degree at most $24$, and
it is permutation
isomorphic to one of the $24$ groups
in Table \ref{badgrouptable}
with their natural permutation
representation, unless indicated otherwise in the table.
\qed

\end{theorem}

For our purposes, we need the following refinement of
Theorem \ref{24special}.
Lemma \ref{factorialpairkeptlem} below is the reason
for the different treatment of degrees
$2,3,4$.

\begin{table}

\centerline{
\begin{tabular}{|l|r|r|r|r|r|}
\hline
Degree $n$	&	Group $G$	&	Order $|G|$	&	$\lg (|G|)/n\leq $	&	$d_n $	\\ 	\hline \hline
5	&	AGL(1,5)	&	20	&	0.8644	&	0.8644	\\	\hline
6	&	PSL(2,5)	&	60	&	0.9845	&		\\	
6	&	PGL(2,5)	&	120	&	1.1512	&	1.1512	\\	\hline
7	&	PSL(3,2)	&	168	&	1.0561	&	1.0561	\\	\hline
8	&	A$\Gamma $L(1,8)	&	168	&	0.9241	&	\\	
8	&	PSL(2,7)	&	168	&	0.9241	&	\\	
8	&	PGL(2,7)	&	336	&	1.0491	&	\\	
8	&	ASL(3,2)=AGL(3,2)	&	1344	&	1.2991	&	1.2991	\\	\hline
9	&	AGL(2,3)	&	432	&	0.9728	&		\\	
9	&	PSL(2,8)	&	504	&	0.9975	&		\\	
9	&	P$\Gamma $L(2,8)	&	1512	&	1.1736	&	1.1736	\\	\hline
10	&	PGL(2,9)	&	720	&	0.9492	&		\\	
10	&	$M_{10}$	&	720	&	0.9492	&		\\	
10	&	$S_6 $ primitive on 10 elt.	&	720	&	0.9492	&	\\	
10	&	P$\Gamma $L(2,9)	&	1440	&	1.0492	&	1.0492	\\	\hline
11	&	$M_{11}$	&	7920	&	1.1774	&	1.1774	\\	\hline
12	&	$M_{11} $ on 12 elements	&	7920	&	1.0793	&		\\	
12	&	$M_{12}$	&	95040	&	1.3781	&	1.3781	\\	\hline
13	&	PSL(3,3)	&	5616	&	0.9582	&	0.9582	\\	\hline
15	&	PSL(4,2)	&	20160	&	0.9533	&	0.9533	\\	\hline
16	&	$2^4 :A_7 $	&	40320	&	0.9563	&	\\	
16	&	$2^4 :L_4 (2) =$AGL(4,2)	&	322560	&	1.1438	&	1.1438	\\	\hline
23	&	$M_{23}$	&	10200960	&	1.0123	&	1.0123	\\	\hline
24	&	$M_{24}$	&	244823040	&	1.1612	&	1.1612	\\	
\hline
\end{tabular}
}

\caption{The primitive permutation groups
given in \cite{Maroti} which have degree $n$, order at least
$2^{n-1} $, and which do not contain the alternating group.
Data from \cite{BuekLee}.
In case of mismatched names, both names
are given.
For each group $G$ that is listed, the upper bound for $\lg (|G|)/n$
is obtained by rounding up to the fourth digit.
For each $n$, we have $d_n \geq \lg (|G|)/n$
for all primitive groups $G$ of degree $n$ that do not contain the
alternating group.
}
\label{badgrouptable}

\end{table}

\begin{define}
\label{defineexponents}

Let $n\in {\mat N} \setminus \{ 1\}
$.
For
$n \in \{ 5,6,7,8,9,10,11,12,13,15,16,23,24\} $,
we set $d_n $ to be equal to the value in
the last column of
Table \ref{badgrouptable}.
For $n\in \{ 2,3,4\} $, we set
$d_2 :={1\over 2} = {\lg (2!)\over 2} $,
$d_3 :=0.8617\geq {\lg (3!)\over 3} $
and
$d_4 :=1.1463\geq {\lg (4!)\over 4} $.
We set
$I
:=\{ 2,3,4;5,6,7,8,9,10,11,12,13,15,16,23,24\} $
and, for $n\not\in I$, we set
$d_n :=2^n $.

\end{define}

\begin{cor}
\label{marotiexpbound}

(Compare with Corollary 1.2 in  \cite{Maroti}, which provides
$|G| < 3^n $.)
Let $G$ be a
primitive subgroup of $S_n $.
If $n\leq 4$,
then
$|G| < 2^{d_n n}$.
If $n\geq 5$
and $G$ does not contain $A_n $,
then
$|G| \leq 2^{d_n n}$.

\end{cor}

{\bf Proof.}
By Theorem \ref{24special},
for $n\not\in I$,
we have that
$|G|\leq 2^{n-1} <2^n =2^{d_n n} $, and,
for $n\in I\setminus \{ 2,3,4\} $,
we have that
$|G|\leq 2^{d_n n} $.
For $n\in \{ 2,3,4\} $,
we have that
$2^{d_n n} \geq n!\geq |G|$.
\qed

\vspace{.1in}

For the following, the author apologizes for
restating some elementary
items from the theory of permutation groups.
However,
standard works on permutation groups, such as \cite{Dixmort},
or references, such as \cite{Maroti},
formulate results
primarily in the language of groups,
which does not seem to have a simple connection to
ordered sets.
The following will
keep the presentation self contained.

\begin{define}

Let $G$ be a permutation group on the set $X$
and let $C\subseteq X$.
We define
$G\cdot C:=\{ \sigma [C]:\sigma \in G\} $.

\end{define}

\begin{prop}
\label{folkloreonblocks}

(Folklore.)
Let $G$ be a
transitive
permutation group on the set $X$
and let $B$ be a $G$-block.
Then
the {\bf block system}
$G\cdot B$ is a partition of $X$,
for every $\sigma \in G$ and every block $B'\in G\cdot B$,
we have
$\sigma [B']\in G\cdot B$, and,
if $A$ is a block that contains $B$, then
every block in $G \cdot B$ that intersects $A$
is contained in $A$.

\end{prop}

{\bf Proof.}
Let $\sigma , \tau \in G$ such that
$\sigma [B]\cap \tau [B]\not= \emptyset $.
Then $\tau ^{-1} \sigma [B]\cap B\not= \emptyset $.
Hence, because $B$ is a block, we infer
$\tau ^{-1} \sigma [B]\cap B=B$, and then $\sigma [B]= \tau [B]$.
This proves that
$G\cdot B$ is a partition of $X$. By definition,
every $\sigma \in G$ maps every block in $G\cdot B$
to another block in $G\cdot B$.

Finally, let $A$ be a block that contains $B$
and let $C:=\sigma [B]\in G\cdot B$ intersect $A$.
Then $\sigma [A]\cap A\not= \emptyset $,
and therefore $\sigma [A]=A$.
Consequently, $C=\sigma [B]\subseteq \sigma [A]=A$.
\qed

\vspace{.1in}

In our investigation of automorphism groups, we will
regularly work with subgroups.
Because a block for a group is also a block for every subgroup,
the following definition is sound.

\begin{define}
\label{actiononblocksys}

Let $G$ be a permutation group on the set $X$
and let $B\subseteq X$ be a block.
For $\sigma \in G$, and $Q\subseteq G\cdot B$,
we define the {\bf action}
$\sigma \hochdamit ^{Q} $ on $Q$ by,
for all $C\in Q$, setting
$\sigma \hochdamit ^{Q} (C):=\sigma [C]$.
Let $A$ be a block that contains $B$
and let $G'$ be a subgroup of $G$.
We define
$A[G\cdot B]:=\{ B'\in G\cdot B: B'\subseteq A\} $
and the {\bf action of $G'$ on $A[G\cdot B]$} as
$G'\hochdamit ^{G\cdot B} _A :=
\{ \sigma \hochdamit ^{G\cdot B} |_{A[G\cdot B]} :
\sigma \in G', \sigma [A]\subseteq A\} $.\footnote{Unless the
author is mistaken, the $G\hochdamit ^{G\cdot B} _A $ are sections of $G$.
However, the definition of sections might include other structures and,
to the author's knowledge, does not explicitly refer to the block structures
we define here.}
When $A$ is equal to the underlying set $X$, we will omit the subscript.

\end{define}

\begin{prop}
\label{folkloreonblocks2}

(Folklore.)
Let $G$ be a
transitive
permutation group on the set $X$
and let $S\subseteq B$ be $G$-blocks.
Then, for any two
$B_1 , B_2 \in G\cdot B$, we have that
$G\hochdamit _{B_1 } ^{G\cdot S} $
is isomorphic to
$G\hochdamit _{B_2 } ^{G\cdot S} $.

\end{prop}

{\bf Proof.}
Let
$B_1 , B_2 \in G\cdot B$
and let
$\mu \in G$ be so that $\mu [B_1 ]=B_2 $.
The function $\Phi (\sigma ):=\mu \sigma \mu ^{-1} |_{B_2 } $
is an isomorphism from
$G\hochdamit _{B_1 } ^{G\cdot S} $
to
$G\hochdamit _{B_2 } ^{G\cdot S} $.
\qed

\vspace{.1in}

Because primitive permutation groups are
the fundamental combinatorial entities from which all
other permutation groups arise,
we define the following.

\begin{define}
\label{primpair}

Let $G$ be a
permutation group
on the set $X$, let
$B\subseteq X$ be a block
and let
$S\subsetneq B$ be a block.
Then $(S,B)$ is called a {\bf primitive
pair} for $G$
iff
$G\hochdamit ^{G\cdot S} _B $
is primitive.
We call ${|B|\over |S|} $ the
{\bf width} of the primitive pair.

\end{define}

\begin{define}

Let $G$ be a permutation group on the set $X$.
Two blocks
$B,B'\subseteq X$ are called
{\bf conjugate blocks}
iff there is a $\sigma \in G$ such that
$B'=\sigma [B]$.
For any primitive pair $(S,B)$ and
any $\sigma \in G $,
we define
$\sigma (S,B):=(\sigma [S],\sigma [B])$.
Two
primitive pairs
$(S,B)$ and $(S',B')$
are called
{\bf conjugate primitive pairs}
iff there is a $\sigma \in G$ such that
$\sigma (S,B)=(S',B')$.
Let
${\cal S}$ be a set of blocks, let ${\cal P}$
be a set of primitive pairs, and let $G'\subseteq G $.
Then we define $G' \cdot {\cal S}
:=
\{ \sigma [S]:S\in {\cal S}, \sigma \in G'\} $
and
$G' \cdot {\cal P}
:=
\{ \sigma (S,B):(S,B)\in {\cal P}, \sigma \in G'\} $.

\end{define}

For a system of conjugate blocks $G\cdot S$, any member can act as
a representative and this will make
estimates of sizes
$\left| G\hochdamit ^{G\cdot S} _B \right| $
easier to read.
Because blocks can be nested within each other,
we define the following.

\begin{define}

Let $G$ be a transitive permutation group
on the set $X$.
A sequence
$\{ x \} =B_0 \subsetneq B_1 \subsetneq
\cdots \subsetneq B_{m} =X$
is called a {\bf primitive nesting}
for $G$ iff, for all $j\in \{ 1, \ldots , m\} $,
we have that $(B_{j-1} , B_j )$ is a primitive pair.

\end{define}

By starting with a singleton $B_0 =\{ x\} $
and successively choosing a smallest size block
that contains, but is not equal to, the block currently under consideration,
we see that every
transitive finite permutation group indeed
has many primitive nestings.
Corollary \ref{marotiexpbound}
shows that, unless
$G\hochdamit ^{G\cdot B_{j-1} } _{B_{j} } $
is an alternating or symmetric group on
$w:={|B_j |\over |B_{j-1} |} \geq 5$
elements,
we have $\left| G\hochdamit ^{G\cdot B_{j-1} } _{B_{j} } \right|
\leq 2^{d_w  w} $.
For groups
$G\hochdamit ^{G\cdot B_j } _{B_{j+1} } $
that contain the alternating group
on their respective domains,
we must devise special treatments
and we define the following.

\begin{define}
\label{factpair}

Let $G$ be a
permutation group
on the set $X$ and let
$(S,B)$ be a primitive
pair for $G$.
Then the primitive pair
$(S,B)$ is called a {\bf factorial
pair} for $G$ iff
$G\hochdamit ^{G\cdot S} _B $
contains
the alternating group $A_{|B[G\cdot S]|} (B[G\cdot S])$
acting
on the set of blocks $B[G\cdot S]$.

\end{define}

\section{Reverse Depth
Enumerations for Primitive Nesting Collections}
\label{controlsetsec}

Absence of nontrivial order-autonomous antichains is
crucial for obtaining exponential upper bounds, and,
for the Automorphism Conjecture,
Proposition \ref{onelexiter}
allows us
to exclude nontrivial order-autonomous antichains
from further consideration.
Moreover, we will later need to
apply some results to subgroups $G^* $ of the
automorphism group,
and we will mostly work with the actions of $G^* $ on orbits.
Therefore, we introduce
the following terminology.

\begin{define}
\label{standardconfdef}

Let $U$
be an ordered set without nontrivial
order-autonomous antichains.
Let $G^* $ be a
subgroup of ${\rm Aut} (U)$.
The set of all $G^* $-orbits will be denoted ${\cal D}$.
We will call $(U,{\cal D} ,G^* )$
a {\bf standard configuration}.

For every $G^* $-orbit
$D\in {\cal D}$,
we define
$\Lambda ^* (D):=
\left\{ \Phi |_D: \Phi \in G^* \right\} $,
where the superscript shall be our indication that we work with
the subgroup $G^* $.

\end{define}

Clearly, the group $G^* $, viewed as a
permutation group,
will not be transitive.
Rather, the groups $\Lambda ^* (D)$
can be viewed as a collection of
``interacting transitive permutation groups."
It will prove useful to
associate a primitive
nesting with every $D\in {\cal D}$.

\begin{define}
\label{primnestcoll}

Let $(U,{\cal D} , G^* )$ be a standard configuration.
For every $D\in {\cal D}$, we can choose a primitive nesting
$\{ x^D \} =B_0 ^D \subsetneq B_1 ^D \subsetneq
\cdots \subsetneq B_{m_D} ^D =D$
for
$\Lambda  ^* (D)$.
Upon making such a choice for every $D\in {\cal D}$,
we call the set
${\cal B} :=\left\{ B_j ^D : D\in {\cal D} , j\in \{ 0,\ldots , m_D \} \right\} $
a {\bf primitive nesting collection}
for $(U,{\cal D},G^* )$.
We let ${\cal F}\subseteq (G^* \cdot {\cal B} )^2 $ be
the set of all factorial pairs $(S,B)\in (G^* \cdot {\cal B} )^2 $
of width $w\geq 5$.

We call
$(U,{\cal D} , G^* ; {\cal B})$, as well as,
when factorial pairs are specifically addressed,
$(U,{\cal D} , G^* ; {\cal B} ,{\cal F} )$, a
{\bf nested standard configuration}.

\end{define}

Within a primitive nesting collection ${\cal B}$,
ordered by inclusion,
for every non-maximal block, we have unique successor.
Because the pair consisting of a
block and its successor is a primitive pair,
these unique successors
are useful for enumeration,
and they
will be useful for all
blocks in $G^* \cdot {\cal B}$, too.

\begin{define}

Let $(U,{\cal D} , G^* ; {\cal B} )$ be a
nested standard configuration.
For every
$C\in G^* \cdot {\cal B} $,
we define $D_C \in {\cal D}$ to be the unique orbit
that
contains $C$.
For every
$C\in {\cal B} $,
we define the
{\bf nesting depth} of $C$ to be the unique index $j$ such that
$C=B_j ^{D_C} $.
Moreover, when the nesting depth of
$C\in {\cal B}$ is
$j<m_{D_C} $, we define
$C^\dag $ to be the unique block in
${\cal B}$
such that $\left( C,C^\dag \right) $ is a primitive pair, that is,
$\left( C,C^\dag \right) =\left( B_j ^{D_C} ,B_{j+1} ^{D_C} \right) $.
For every conjugate $\widetilde{C} \in
\Lambda  ^* (D_C  ) \cdot C$
of $C$, we define the {\bf successor}
$\widetilde{C} ^\dag \in
\Lambda  ^* (D_C  ) \cdot C^\dag $
to be the unique conjugate of $C^\dag $
such that $\left( \widetilde{C} , \widetilde{C} ^\dag \right) $
is a primitive pair.

\end{define}

Some of the following
lemmas will be needed later for general sets ${\cal K}$ of
$G^* $-blocks that need not be
associated with a primitive nesting ${\cal B}$.
Hence, although in this section these sets will
satisfy ${\cal K} \subseteq {\cal B}$,
when possible, results will be proved in the absence of
any
primitive nesting collections.

\begin{define}
\label{restrictonK}

Let $(U,{\cal D} , G^*
)$ be a
standard configuration,
let $G'$ be a subgroup of $G^* $,
and let
${\cal K}
$ be a set of blocks for $G^* $.
We define
$$G' _{\setcontrol {\cal K} }
:=\left\{
\Phi \in G' :
(\forall K\in {\cal K} )
\Phi \hochdamit ^{\Lambda  ^* (D_K )\cdot K }
=
{\rm id} \hochdamit ^{\Lambda  ^* (D_K )\cdot K }
\right\} .
$$
In cases in which ${\cal K}$ can be easily enumerated
as ${\cal K} =\{ B_1 , \ldots , B_n \} $, we will
write
$G' _{\setcontrol B_1 , \ldots , B_n } $
instead of
$G' _{\setcontrol \{ B_1 , \ldots , B_n \} } $.

\end{define}

\begin{lem}
\label{normalgrouplem}

Let $(U,{\cal D} , G^*
)$ be a
standard configuration,
let $G'$ be a subgroup of $G^* $,
and let
${\cal K}
$ be a set of blocks for $G^* $.
Then $G' _{\setcontrol {\cal K} } $
is a normal subgroup of $G' $.
Moreover,
for every
$B\not\in
{\cal K} $, the factor
group
$G' _{\setcontrol {\cal K} } /
G' _{\setcontrol {\cal K} \cup \{ B\} } $
is isomorphic to
$\left\{
\Phi
\hochdamit ^{\Lambda  ^* (D_B )\cdot B }
:
\Phi \in
G' _{\setcontrol {\cal K} }
\right\}
$.

\end{lem}

{\bf Proof.}
Let $\Phi \in G' _{\setcontrol {\cal K} } $,
let $\Psi \in G' $ and let
$K\in {\cal K}$. Then
\begin{eqnarray*}
\left( \Psi ^{-1}
\circ
\Phi
\circ
\Psi \right)
\hochdamit ^{
\Lambda  ^* (D_K )\cdot K  }
& = &
\left( \Psi ^{-1} \right)
\hochdamit ^{
\Lambda  ^* (D_K )\cdot K  }
\circ
\Phi
\hochdamit ^{
\Lambda  ^* (D_K )\cdot K  }
\circ
\Psi
\hochdamit ^{
\Lambda  ^* (D_K )\cdot K }
\\
& = &
\left( \Psi ^{-1} \right)
\hochdamit ^{
\Lambda  ^* (D_K )\cdot K }
\circ
{\rm id}
\hochdamit ^{
\Lambda  ^* (D_K )\cdot K }
\circ
\Psi
\hochdamit ^{
\Lambda  ^* (D_K )\cdot K }
\\
& = &
\left( \Psi ^{-1} \right)
\hochdamit ^{
\Lambda  ^* (D_K )\cdot K }
\circ
\Psi
\hochdamit ^{
\Lambda  ^* (D_K )\cdot K }
\\
& = &
\left( \Psi ^{-1} \Psi \right)
\hochdamit ^{
\Lambda  ^* (D_K )\cdot K  }
\\
& = &
{\rm id}
\hochdamit ^{
\Lambda  ^* (D_K )\cdot K  } ,
\end{eqnarray*}
and hence, because $K\in {\cal K}$ was arbitrary,
$G' _{\setcontrol {\cal K} } $
is a normal subgroup of
$G' $.

Now let
$B
\not\in
{\cal K} $
and let
$[\Phi ], [\Psi ]\in G' _{\setcontrol {\cal K} } /
G' _{\setcontrol {\cal K} \cup \{ B\} } $.
Then $[\Phi ]=[\Psi ]$
iff there is a $\Delta \in G' _{\setcontrol {\cal K} \cup \{ B\} } $
such that
$\Phi
\circ
\Delta =\Psi $,
which is the case iff
$
\Phi
\hochdamit ^{\Lambda  ^* (D_B )\cdot B }
=
\Psi
\hochdamit ^{\Lambda  ^* (D_B )\cdot B }
$.
Thus the function
$
[\Phi ]\mapsto
\Phi
\hochdamit ^{\Lambda  ^* (D_B )\cdot B }
$
is well-defined and injective
from
$G' _{\setcontrol {\cal K} } /
G' _{\setcontrol {\cal K} \cup \{ B\} } $
to
$\left\{
\Phi
\hochdamit ^{\Lambda  ^* (D_B )\cdot B }
:
\Phi \in G' _{\setcontrol {\cal K} }
\right\}
$.
Because the function
clearly is surjective and
a homomorphism, we obtain that
$G' _{\setcontrol {\cal K} } /
G' _{\setcontrol {\cal K} \cup \{ B\} } $
is isomorphic to
$\left\{
\Phi
\hochdamit ^{\Lambda  ^* (D_B )\cdot B }
:
\Phi \in
G' _{\setcontrol {\cal K} }
\right\}
$.
\qed

\begin{define}
\label{fixedbyprevious}

Let $(U,{\cal D} , G^* ; {\cal B})$ be a nested
standard configuration,
let
${\cal K} \subseteq {\cal B}$
and let $B\in {\cal B}\setminus {\cal K}$.
We define
\begin{eqnarray*}
\Lambda  ^* \left( {\cal K}\not\setcontrol B \right)
& := &
\left\{
\Phi
\hochdamit ^{\Lambda  ^* (D_B )\cdot B }
:
\Phi \in
G' _{\setcontrol {\cal K} }
\right\}
.
\end{eqnarray*}

\end{define}

\begin{remark}
\label{simpleupperbound}

{\rm
Because
$
\left|
\Lambda  ^* \left( {\cal K}\not\setcontrol B \right)
\right|
$
can be cumbersome to compute,
via Proposition \ref{folkloreonblocks2},
we
note the following upper bounds when $B^\dag \in {\cal K}$.
\begin{eqnarray*}
\left|
\Lambda  ^* \left( {\cal K}\not\setcontrol B \right)
\right|
& \leq &
\left|
\left\{
\Phi
\hochdamit ^{\Lambda  ^* (D_B )\cdot B }
:
\Phi \in G^* ,
\Phi
\hochdamit ^{\Lambda  ^* (D_B )\cdot B^\dag }
=
{\rm id}
\hochdamit ^{\Lambda  ^* (D_B )\cdot B^\dag }
\right\}
\right|
\\
& \leq &
\left|
\Lambda  \hochdamit ^{\Lambda  ^* (D_B )\cdot B }
_{B^\dag }
\right| ^{
\left|
D_B
\right|
\over
\left|
B^\dag
\right|
}
\end{eqnarray*}
Note that, when
$\left( B, B^\dag \right) $
is not a factorial pair of width $\geq 5$, then the above is
bounded by
$\left(
2^{d_{\left|
B^\dag
\right|
/ |B |}
{
\left|
B^\dag
\right|
\over
\left|
B
\right|
}
}
\right)
^{
\left|
D_B
\right|
\over
\left|
B^\dag
\right|
}
=
2^{d_{\left|
B^\dag
\right|
/ |B |}
{
\left|
D_B
\right|
\over
\left|
B
\right|
}
}
$.
Products of such powers will be crucial for our estimates and they
lead us to Lemma \ref{1.54bound} below.
\qex
}

\end{remark}

\begin{lem}
\label{1.54bound}

Let $1\leq c_0 <c_1 <\cdots <c_m =x$ be natural numbers such that,
for all $j\in \{ 0,\ldots , m-1\} $, the quotient
${c_{j+1} \over c_j } $
is an integer greater than or equal to $2$.
Then
$
\prod _{j=0} ^{m-1}
2^{ d_{
c_{j+1}
/ c_j }
{
x
\over
c_j }
}
\leq
2^{1.54 {x\over c_0 } } .$

\end{lem}

{\bf Proof.}
Note that, for $b_j :={c_j \over c_0 } $, we have that
${b_{j+1} \over b_j }= {c_{j+1} \over c_j } \geq 2$.
Because $b_0 =1$, we obtain $b_j \geq 2^j $.
Now we have the following.
\begin{eqnarray*}
\lg
\left(
\prod _{j=0} ^{m-1}
2^{ d_{
c_{j+1}
/ c_j }
{
x
\over
c_j }
}
\right)
& = &
\sum _{j=0} ^{m-1}
d_{
c_{j+1}
/ c_j }
{
x
\over
c_j }
=
{x\over c_0 }
\sum _{j=0} ^{m-1}
d_{
b_{j+1}
/ b_j }
{1\over b_j }
\\
& \leq &
{x\over c_0 }
\left(
\sum _{j=0} ^s
d_{
b_{j+1}
/ b_j }
{1\over b_j }
+
\sum _{j=s+1} ^{m-1}
{
d_{
b_{j+1}
/ b_j }
\over b_j }
\right)
\\
& \leq &
{x\over c_0 }
\left(
\sum _{j=0} ^s
d_{
b_{j+1}
/ b_j }
{1\over b_j }
+
\sum _{j=s+1} ^{m-1}
{1.38
\over 2^j}
\right)
\\
& \leq &
{x\over c_0 }
\left(
\sum _{j=0} ^s
d_{
b_{j+1}
/ b_j }
\prod _{i=1} ^j {b_{i-1} \over b_i }
+
{1.38
\over 2^s}
\right)
\end{eqnarray*}

It is simple to write the nested loops that compute
an upper bound for
$\sum _{j=0} ^s
d_{
b_{j+1}
/ b_j }
\prod _{i=1} ^j {b_{i-1} \over b_i }
$
for $s=7$:
The factor
$d_{b_{j+1}
/ b_j } $
from Definition \ref{defineexponents}
only depends on
the number
$
b_{j+1}
/ b_j $,
and, for
$
b_{j+1}
/ b_j \not\in I$, we have
$d_{
b_{j+1}
/ b_j  }
=1$.
For any reciprocal
${b_i \over b_{i-1} }
\in
I$,
we use
${b_{i-1} \over b_i }$ as a factor in the product.
For any reciprocal
${b_i \over b_{i-1} }
\not\in
I
$,
we use
${1\over 14} $ as an upper bound
for the factor
${b_{i-1} \over b_i }$
in the product.
Thus, for each factor,
if the
reciprocal is not in $I$,
we can use a single number,
leading to
finitely many combinations that need to be considered.

Computationally, we obtain
$\sum _{j=0} ^7
d_{
b_{j+1}
/ b_j }
\prod _{i=1} ^j {b_{i-1} \over b_i }
\leq
1.5284$.
Together with
${1.38\over 2^7 } \leq
0.0108$,
we obtain
${x\over c_0 }
( 1.5284+0.0108) \leq 1.54 {x\over c_0 }  $
as an upper bound for the logarithm.
\qed

\vspace{.1in}

Lemma
\ref{normalgrouplem}
suggests the successive expansion of
a set ${\cal K}\subset {\cal B}$ by adding controlled blocks.
To apply Lemma \ref{1.54bound}, we must have primitive
pairs, which is the motivation behind
part \ref{RDEdef2} of Definition \ref{RDEdef}
below.

\begin{define}
\label{RDEdef}

Let $(U,{\cal D} , G^* ; {\cal B})$ be a nested
standard configuration.
A vector $\overrightarrow{{\cal K}} =(B_1 , \ldots , B_k )\in {\cal B} ^k $
is a {\bf reverse depth enumeration}
iff
the following hold.
\begin{enumerate}
\item
\label{RDEdef1}
For all distinct $i, j\in \{ 1, \ldots , k\} $,
we have $B_i \not= B_j $.
\item
\label{RDEdef2}
For all
$j\in \{ 1, \ldots , k\} $, the set
${\cal K} _j :=\{ B_1 , \ldots , B_j \} $
is {\bf upper block closed}, that is,
for all $K\in {\cal K}_j $ and $B\in {\cal B}$, we have that
$K\subseteq B$ implies $B\in {\cal K}_j $.
\end{enumerate}
We set
${\cal K} :=\{ B_1 , \ldots , B_k \} $, and,
for every
$B \in {\cal K} $,
we define
$\overleftarrow{\cal K} ( B ):=
{\cal K}_{j-1 } $, where
$j$ is the index of $B$ in the reverse depth enumeration.
For simpler notation, we will write
$
\Lambda  ^* \left( \overleftarrow{\cal K} \not\setcontrol B \right)
$
instead of
$
\Lambda  ^* \left( \overleftarrow{\cal K} (B)\not\setcontrol B \right)
$.

We call
$\left( U,{\cal D} , G^* ; {\cal B}, \overrightarrow{{\cal K}} \right)$
a {\bf partially enumerated nested
standard configuration}.
In case ${\cal K}={\cal B}$, we will
call
$\left( U,{\cal D} , G^* ; \overrightarrow{\cal B}\right)$,
as well as
$\left( U,{\cal D} , G^* ; \overrightarrow{\cal B}, {\cal F}\right)$,
an
{\bf enumerated nested
standard configuration}.

\end{define}

Components of a reverse depth enumeration $\overrightarrow{K} $ will
typically be denoted $C$, as they are
considered to be {\em controlled}.

\begin{prop}
\label{controlbysubsetwithgroups}

Let
$\left( U,{\cal D} , G^* ; \overrightarrow{{\cal B}} \right)$
be an enumerated nested
standard configuration.
Then
$\displaystyle{
\left| G^* \right|
=
\prod _{C\in {\cal B} }
\left|
\Lambda  ^* \left( \overleftarrow{\cal B} \not\setcontrol C \right)
\right|
.
} $

\end{prop}

{\bf Proof.}
Let
$\overrightarrow{{\cal B}} =:(B_1 , \ldots , B_m )$
and let $B_0 :=\emptyset $.
%and, for every $k\in \{0, \ldots , m\} $,
%define
%${\cal B}_k  :=\{ B_1 , \ldots , B_k \} $.
By Lemma \ref{normalgrouplem},
for every
$k\in \{1, \ldots , m\} $, the factor
group
$G' _{\setcontrol {\cal B}_{k-1} } /
G' _{\setcontrol {\cal B}_k }
=
G' _{\setcontrol {\cal B}_{k-1} } /
G' _{\setcontrol {\cal B}_{k-1} \cup \{ B_k \} }
$
is isomorphic to
$
\Lambda  ^* \left( {\cal B}_{k-1} \not\setcontrol B_k \right)
=
\Lambda  ^* \left( \overleftarrow{\cal B} \not\setcontrol B_k \right)
$.

Because $\{ {\rm id} \} =
G^* _{\setcontrol {\cal B} _m }
\triangleleft G^* _{\setcontrol {\cal B} _{m-1} }
\triangleleft
\cdots
\triangleleft G^* _{\setcontrol {\cal B} _{1} }
\triangleleft G^* _{\setcontrol {\cal B} _{0} }
=G^* $
is a normal series, we obtain
$\left| G^* \right|
=
\prod _{k=1} ^m
\left|
G' _{\setcontrol {\cal B}_{k-1} } /
G' _{\setcontrol {\cal B}_k }
\right|
=
\prod _{k=1} ^m
\left|
\Lambda  ^* \left( \overleftarrow{\cal B} \not\setcontrol B_k \right)
\right|
$.
\qed

\vspace{.1in}

If none of the primitive groups
$\Lambda  ^* \left( \overleftarrow{\cal B} \not\setcontrol C \right) $
contains an alternating group, then
Remark \ref{simpleupperbound} and
Lemma \ref{1.54bound}
provide an exponential bound
for the product in
Proposition
\ref{controlbysubsetwithgroups}, albeit with an exponent that is
slightly too high.
To gain better control over the product, we
introduce a bounding product in which
every factor
$\gamma _{\cal K} ^* (C)$
corresponds to a
$
\left|
\Lambda  ^* \left( \overleftarrow{\cal B} \not\setcontrol C \right)
\right|
$, but in which we
will have the freedom
to
``trade factors
between sets of blocks."
There will be three types of
problematic factors $\gamma _{\cal K} ^* (C)$,
and they give rise to the set ${\cal K}_u ^\gamma $ below:
Those that do not satisfy a
bound like a primitive group of the same degree would satisfy
and
those for which the block $C$ is of level zero, but
insufficient to
drop the factor in the exponent from
$1.54$ to $0.99$, in which case we distinguish
whether
$\left( C, C^\dag \right) $ is factorial of
width $\geq 5$ or not.

\begin{define}
\label{controlfunction}

Let
$\left( U,{\cal D} , G^* ; {\cal B}, \overrightarrow{{\cal K}} \right)$
be a partially enumerated nested
standard configuration.
We call the function
$\gamma _{\cal K} ^* :{\cal K} \to [1, \infty )$
a {\bf control function}
iff
the following hold.
\begin{enumerate}
\item
\label{controlfunction1}
$\displaystyle{
\prod _{C\in {\cal K} }
\left|
\Lambda  ^* \left( \overleftarrow{\cal K} \not\setcontrol C \right)
\right|
\leq
\prod _{C\in {\cal K} }
\gamma _{\cal K} ^* (C)
.}$
\item
\label{controlfunction2}
For all
$C\in  {\cal K} $ of
depth at least $1$
such that
$\left( C, C^\dag \right) $ is
of width $w$, but
not factorial
of width $w\geq 5$,
we have
$
\gamma _{\cal K} ^* (C)
\leq
\left|
\Lambda  ^* \left( \overleftarrow{\cal K} \not\setcontrol C \right)
\right|
\leq
2^{d_w
|\Lambda  ^* (D_C )\cdot C|} $.
\end{enumerate}
We also define
${\cal K}_f $
to be the set of all
$S\in  {\cal K} $ such that
$\left( S, S^\dag \right) $ is factorial
of width $w\geq 5$,
we define
${\cal K}_0 $
to be the set of all
$S\in  {\cal K} $ that
have depth $0$,
and we define the set of
{\bf ($
\gamma _{\cal K} ^*
$-)uncontrolled factors in ${\cal K}$} to be
\begin{eqnarray*}
{\cal K}_u ^\gamma
& := &
\left\{ C\in {\cal K}_f \setminus {\cal K} _0 :
\gamma _{\cal K} ^* (C)
>
2^{d_{\left| C^\dag \right| /|C|}
|\Lambda  ^* (D_C )\cdot C|}
\right\}
\\
& &
\cup
\left\{ C\in {\cal K} _f \cap {\cal K}_0 :
\gamma _{\cal K} ^* (C)
\prod _{T\subset D_C , T\in {\cal K} \setminus {\cal K} _f }
\gamma _{\cal K} ^* (T)
>
2^{
0.99
|D_C | }
\right\}
\\
& &
\cup
\left\{ C\in {\cal K} _0 \setminus {\cal K}_f :
\prod _{T\subset D_C , T\in {\cal K} \setminus {\cal K} _f }
\gamma _{\cal K} ^* (T)
>
2^{
0.99
|D_C | }
\right\}
.
\end{eqnarray*}

\end{define}

Because our combinatorial analysis will focus on
individual orbits and their interactions, the
following reindexing will be helpful.

\begin{define}
\label{controlforsubset}

Let
$\left( U,{\cal D} , G^* ; {\cal B}, \overrightarrow{{\cal K}} \right)$
be a partially enumerated nested
standard configuration,
and let $\gamma _{\cal K} ^* $ be a control function.
For every $D\in {\cal D}$, we define
${\cal K} ^D := {\cal K} \cap \left\{ B_0 ^D , B_1 ^D ,
\ldots , B_{m_D -1} ^D \right\} $,
${\cal K} ^{\gamma ,D} _u := {\cal K} _u ^\gamma
\cap \left\{ B_0 ^D , B_1 ^D ,
\ldots , B_{m_D -1} ^D \right\} $,
\begin{eqnarray*}
\Omega _{{\cal K} } \left( D\right)
:=
\prod _{C\in {\cal K} ^D }
\gamma _{\cal K} ^* (C)
& {\rm and} &
\Upsilon _{{\cal K} } \left( D\right)
:=
\prod _{C\in {\cal K} ^{\gamma ,D} _u }
\gamma _{\cal K} ^* (C)
.
\end{eqnarray*}

\end{define}

\begin{prop}
\label{reindexed}

Let
$\left( U,{\cal D} , G^* ; \overrightarrow{{\cal B}} \right)$
be an enumerated nested
standard configuration,
and let $\gamma _{\cal B} ^* $ be a control function.
Then
$\displaystyle{
\left| G^* \right|
\leq
\prod _{D\in {\cal D}}
\Omega _{{\cal B} } \left( D\right)
\leq
2^{0.99|U|}
\prod _{D\in {\cal D}}
\Upsilon _{{\cal B} } \left( D\right)
.
}$

\end{prop}

{\bf Proof.}
By Proposition \ref{controlbysubsetwithgroups}
and by part \ref{controlfunction1}
of Definition \ref{controlfunction},
we have
$\left| G^* \right|
=
\prod _{C\in {\cal B} }
\left|
\Lambda  ^* \left( \overleftarrow{\cal B} \not\setcontrol C \right)
\right|
\leq
\prod _{C\in {\cal K} }
\gamma _{\cal B} ^* (C)
$.
The inequalities follow because $\prod _{D\in {\cal D}}
\Omega _{{\cal B} } \left( D\right) $
is a
simple reindexing of
$\prod _{C\in {\cal B} }
\gamma _{\cal B} ^* (C)
$
and from the fact that
$\Omega _{\cal B} (D)\leq 2^{0.99|D|} \Upsilon _{\cal B} (D)$.
\qed

\vspace{.1in}

Finally, we note that the
net contribution of the factors that
are not in ${\cal K} _u ^{\gamma , D} $
is exponential, albeit with a factor that is slightly
too large.

\begin{prop}
\label{level0elimlem}

Let
$\left( U,{\cal D} , G^* ; {\cal B}, \overrightarrow{{\cal K}} \right)$
be a partially enumerated nested
standard configuration,
let $\gamma _{\cal K} ^* $ be a control function,
and let $D\in {\cal D}$.
\begin{enumerate}
\item
\label{level0elimlem1}
If
$B_0 ^D \not\in {\cal K} ^D $,
then
$\Omega _{{\cal K} } (D)\leq
2^{ 1.54 {|D|\over \left| B_1 ^D \right| } }
\Upsilon _{{\cal K} } \left( D\right)
$.
\item
\label{level0elimlem2}
If
$B_0 ^D \not\in {\cal K} ^D $ and
$\left|
\Lambda  ^* \left( {\cal K}
\not\setcontrol B_0 ^D \right)
\right|
\leq
2^{d_{\left| B_{1} ^D \right| /
\left| B_{0} ^D \right| } |D|}
$,
then
\\
$\Omega _{{\cal K} } (D)
\left|
\Lambda  ^* \left( {\cal K}
\not\setcontrol B_0 ^D \right)
\right|
\leq
2^{ 1.54 {|D|} }
\Upsilon _{{\cal K} } \left( D\right)
$.
\end{enumerate}

\end{prop}

{\bf Proof.}
Let ${\cal K} ^D =:\{ C_0 , \ldots , C_{m-1} \} $, enumerated
so that $C_0 \subset C_1 \subset \cdots \subset C_{m-1} \subset C_m :=D $.

To prove part \ref{level0elimlem1},
let $C_0 =B_0 ^D \not\in {\cal K} ^D $.
By Definition
\ref{controlforsubset}
and by Lemma \ref{1.54bound},
we have
\begin{eqnarray*}
{\Omega _{{\cal K} } \left( D\right)
\over
\Upsilon _{{\cal K} } \left( D\right)
}
& = &
\prod _{C\in {\cal K} ^D \setminus
{\cal K} ^{\gamma ,D} _u
}
\gamma _{\cal K} ^* (C)
\leq
\prod _{C\in
\{ C_1 , \ldots , C_{m-1} \}
}
2^{d_{\left|
C^\dag
\right|
/ |C |}
|\Lambda  ^* (D )\cdot C|
}
\\
& \leq &
\prod _{j=1} ^{m-1}
2^{d_{\left|
C_{j+1}
\right|
/ |C_j |}
{|D|\over |C_{j} |}
}
\leq
2^{1.54 {|D|\over |C_1 |} }
\end{eqnarray*}

To prove part \ref{level0elimlem2},
let
$B_0 ^D \not\in {\cal K} ^D $ and
$\left|
\Lambda  ^* \left( {\cal K}
\not\setcontrol B_0 ^D \right)
\right|
\leq
2^{d_{\left| B_{1} ^D \right| /
\left| B_{0} ^D \right| } |D|}
$.
The same computation as for
part \ref{level0elimlem1},
with estimates starting with
$C_0 $ instead of $C_1 $ yields
$\Omega _{{\cal K} } (D)
\left|
\Lambda  ^* \left( {\cal K}
\not\setcontrol B_0 ^D \right)
\right|
\leq
2^{ 1.54 {|D|} }
\Upsilon _{{\cal K} } \left( D\right)
$.
\qed

\vspace{.1in}

For a
fully enumerated nested
standard configuration
$\left( U,{\cal D} , G^* ; \overrightarrow{{\cal B}} \right)$,
Propositions \ref{reindexed}
and \ref{level0elimlem}
show that a control function
$\gamma _{\cal B} ^* $
such that
$\Upsilon _{{\cal B} } \left( D\right) =1$
for all $D\in {\cal D}$
will give the desired bound
$
\left| G^*
\right|
\leq
\prod _{D\in {\cal D}}
\Omega _{{\cal B} } \left( D\right)
\leq
2^{0.99|U|} $.
Note that, for {\em any} reverse depth enumeration $\overrightarrow{\cal K}$ that
does not contain any blocks of depth $0$,
$\gamma _{\cal K} ^* (C):=
\left|
\Lambda  ^* \left( \overleftarrow{\cal K} \not\setcontrol C \right)
\right|
$ defines a
control function.
Therefore,
to achieve a bound
$|{\rm Aut}  (U)|
\leq
2^{0.99|U|} $, we must compensate for
the factorial factors
$\left|
\Lambda  ^* \left( \overleftarrow{\cal B} \not\setcontrol C \right)
\right|
$
that are collected in the product
$\prod _{D\in {\cal D}} \Upsilon _{{\cal B} } \left( D\right) $,
and
we must expand the function to blocks at depth $0$
so that any new large factors can be similarly compensated.
Before we turn to the combinatorial side of this issue, we
record a useful fact about the dichotomous heredity of
factorial pairs.
Because the ideas of
Lemma
\ref{factorialpairkeptlem}
below
will need to be iterated, instead of $G^* $,
we use
a subgroup $G'$ of $G^* $.

\begin{lem}
\label{factorialpairkeptlem}

Let $(U,{\cal D} , G^*
)$ be a
standard configuration,
let $G'$ be a subgroup of $G^* $,
let
${\cal K}
$ be a set of blocks for $G^* $,
and let $(S , B )$
be a factorial
pair of width $w\geq 5$ for $G'$.
If
$G_{\setcontrol {\cal K} } '
\hochdamit ^{\Lambda  ^* (D)\cdot S } $
has more than one element, then
$\left( S ,B \right) $ is a factorial
pair of width $w\geq 5$ for
$G_{\setcontrol {\cal K} } '
$,
too.

\end{lem}

{\bf Proof.}
Because
$\left( S ,B \right) $ is a factorial
pair
of width $w\geq 5$
for $G'$, the group
$G'
\hochdamit ^
{\Lambda  ^* (D )\cdot S } _{B } $
is the symmetric or the alternating group on
the set
$
B [\Lambda  ^* (D )\cdot S ]$, which has
$w\geq 5$ elements.
By Proposition \ref{folkloreonblocks2},
because
$G_{\setcontrol {\cal K} } '
\hochdamit ^{\Lambda  ^* (D)\cdot S } $
has more than one element,
$G_{\setcontrol {\cal K} } '
\hochdamit ^{\Lambda  ^* (D )\cdot S } _{B } $
must have more than one element, too.

We now prove that
$G_{\setcontrol {\cal K} } '
\hochdamit ^{\Lambda  ^* (D )\cdot S } _{B } $
is a normal subgroup of
$G'
\hochdamit ^
{\Lambda  ^* (D )\cdot S } _{B } $.
Let $\sigma
\in G_{\setcontrol {\cal K} } '
\hochdamit ^{\Lambda  ^* (D )\cdot S } _{B } $
and let
$\mu \in
G'
\hochdamit ^
{\Lambda  ^* (D )\cdot S } _{B} $.
Then there are
$\Phi \in
G_{\setcontrol {\cal K} } ' $
such that
$\sigma =
\Phi \hochdamit ^{\Lambda  ^* (D )\cdot S } _{B } $
and
$\Psi \in G' $ such that
$\mu =
\Psi \hochdamit ^{\Lambda  ^* (D )\cdot S } _{B } $.
Now
\begin{eqnarray*}
\mu ^{-1}
\circ
\sigma
\circ
\mu
& = &
\left( \Psi
\hochdamit ^{\Lambda  ^* (D )\cdot S } _{B }
\right) ^{-1}
\circ
\Phi
\hochdamit ^{\Lambda  ^* (D )\cdot S } _{B }
\circ
\Psi
\hochdamit ^{\Lambda  ^* (D )\cdot S } _{B }
\\
& = &
\left( \Psi ^{-1}
\right)
\hochdamit ^{\Lambda  ^* (D )\cdot S } _{B }
\circ
\Phi
\hochdamit ^{\Lambda  ^* (D )\cdot S } _{B }
\circ
\Psi
\hochdamit ^{\Lambda  ^* (D )\cdot S } _{B }
\\
& = &
\left( \Psi ^{-1}
\circ
\Phi
\circ
\Psi
\right)
\hochdamit ^{\Lambda  ^* (D )\cdot S } _{B }
\\
& \in &
G_{\setcontrol {\cal K} } ^*
\hochdamit ^{\Lambda  ^* (D )\cdot S } _{B } ,
\end{eqnarray*}
where the final step is guaranteed by Lemma \ref{normalgrouplem}.
Hence
$G_{\setcontrol {\cal K} } '
\hochdamit ^{\Lambda  ^* (D )\cdot S } _{B } $
is a normal subgroup of
$G'
\hochdamit ^
{\Lambda  ^* (D )\cdot S } _{B } $.

The only normal subgroups of a symmetric group
$S_w $ with $w\geq 5$ are $S_w $ itself,
the alternating group $A_w $, and the
singleton group $\{ {\rm id} \} $.
The only normal subgroups of an alternating group
$A_w $ with $w\geq 5$ are $A_w $ itself and the
singleton group $\{ {\rm id} \} $.
Because
$G_{\setcontrol {\cal K} } '
\hochdamit ^{\Lambda  ^* (D )\cdot S } _{B } $
is not the singleton group $\{ {\rm id} \} $,
it must be the alternating group or the symmetric group
on $
B [\Lambda  ^* (D )\cdot S ]$.
Hence
$\left( S ,B \right) $ is a factorial
pair for
$G_{\setcontrol {\cal K} } '
$.
\qed

\section{Inter-Orbit Combinatorial Control}
\label{combcontsec}

For two orbits $D,E\in {\cal D}$,
the action of
$G^*$
on $E$ can partially or completely
determine its action on $D$.
Lemma \ref{controllemma} and
Definition \ref{controldef}
below
show
how this
``direct inter-orbit combinatorial control"
actually
manifests between block
systems
of the groups
$\Lambda  ^* (E)$
and
$\Lambda  ^* (D)$.
This ``control"
will be our primary tool to avoid, or
compensate for,
large factors
$\gamma _{\cal K} ^* (C)$
that are collected
in $\prod _{D\in {\cal D}} \Upsilon _{\cal K} (D)$.
Because the interactions investigated in this section are local, we
will not fix a global primitive nesting collection.
Rather, a (controlling) block $C\subseteq E$ will
be used to control the
action on (sub)blocks $S \subset D$, which are themselves
contained in a block $B\subseteq D$.
As in the previous section,
blocks will routinely be used as representatives
of their corresponding block systems.
We start by indicating how blocks can be related to each other.

\begin{define}

Let $(U,{\cal D} , G^* )$ be a standard configuration,
let $D,E\in {\cal D}$,
let $S $ be a block of
$\Lambda  ^* (D)$,
and
let
$C$ be a block of
$\Lambda  ^* (E)$.
\begin{enumerate}
\item
We say that $C$ is
{\bf
(nontrivially) woven
with} $S$
and write $C\curlyvee S$ iff there
are elements $c\in C$ and $s\in S$ that are comparable, as well as
elements $c'\in C$ and $s'\in S$ that are incomparable.
\item
We say that
$S$ is {\bf made} by $C$
and write $C\diamondsuit S$
iff $C>S$ or $C<S$.
\item
We say that $C$ is {\bf incomparable to (every element in)} $S$
and write $C\| S$ iff
no element $c\in C$ is comparable to any element $s\in S$.
\end{enumerate}
In case
$S\diamondsuit C$ or $S\| C$, we will also say that
$C$ is {\bf (trivially) woven} with $S$.

\end{define}

Some
weavings
can be considered similar/equivalent.

\begin{define}

Let $(U,{\cal D} , G^* )$ be a standard configuration,
let $D,E\in {\cal D}$, let $S_1 ,S_2 $ be blocks of
$\Lambda  ^* (D)$,
and
let
$C_1 , C_2 $ be blocks of
$\Lambda  ^* (E)$.
Then we say that $(C_1 , S_1 )$ and $(C_2 ,S_2 )$ are
{\bf similarly
woven}
and write $(C_1 , S_1 )\bowtie (C_2 , S_2 )$
iff there is a
$\Phi \in G^*$ such that
$\Phi [C_1 ]=C_2 $ and $\Phi [S_1 ]=S_2 $.

\end{define}

\begin{lem}

Let $(U,{\cal D} , G^* )$ be a standard configuration,
let
$D,E\in {\cal D}$, let
$S$ be a
block of
$\Lambda  ^* (D)$,
and
let
$C$ be a block of
$\Lambda  ^* (E)$.
Then the relation
$\bowtie $ is an equivalence relation on
the product
$\Lambda  ^* (D)\cdot S\times \Lambda  ^* (E)\cdot C$.
Moreover,
in case $S\not= D$ or $C\not= E$,
no $\bowtie $-equivalence class is a singleton.

\end{lem}

{\bf Proof.}
Because
the identity is an automorphism,
$\bowtie $ is reflexive.
Because
every
automorphisms has an
inverse,
$\bowtie $ is symmetric.
Because
the composition of two automorphisms is again an automorphism,
$\bowtie $ is transitive.

Finally, in case $S\not= D$ or $C\not= E$,
one of
$\Lambda  ^* (D)\cdot S$
or
$\Lambda  ^* (E)\cdot C$
has more than one element.
Because $G^* $ acts transitively on
$\Lambda  ^* (D)\cdot S$
and on $\Lambda  ^* (E)\cdot C$,
we infer that
no $\bowtie $-equivalence class is a singleton.
\qed

\begin{define}

Let $(U,{\cal D} , G^* )$ be a standard configuration,
let $D,E\in {\cal D}$, let
$S$ be a
block of
$\Lambda  ^* (D)$,
and
let
$C$ be a block of
$\Lambda  ^* (E)$.
Let ${\cal S}_1 , \ldots , {\cal S} _p $ be the
equivalence classes of
$\bowtie $
on $\Lambda  ^* (D)\cdot S\times \Lambda  ^* (E)\cdot C$
that contain
nontrivially
woven
pairs
$C'\curlyvee S'$.
Then,
for all $j\in \{ 1, \ldots , p\} $ and all
$(C',S')\in {\cal S} _j $,
we will write $C'\curlyvee _j S'$.
These relations technically are
not symmetric, and they, as well as their number $p$,
depend on the
containing orbits, as well as the
block systems involved.
Because the direction, the containing
orbits and the block systems
will be clear from the context, similar to
$\diamondsuit $ and $\| $, they
will not be indicated in the
relation symbol
$\curlyvee _j $
or in the number $p$.

\end{define}

\begin{nota}

{\rm From now on, $p$ will denote the
number of $\bowtie $-equivalence classes that contain
similarly
nontrivially woven
pairs
$C'\curlyvee S'$, and
$\curlyvee _1 , \ldots , \curlyvee _p $
will denote the corresponding weaving relations.
}

\end{nota}

Now note that
the
blocks in one orbit, which have the same
weavings
with
certain blocks in
other
orbits,
can be combined
to form a ``block of blocks."

\begin{define}

Let $w,m\in {\mat N}$ and let $S_1, \ldots , S_w $ be
pairwise disjoint nonempty sets.
We define $(P_1 , \ldots , P_m )$ to be a
{\bf sorted partition}
of $\{ S_1 , \ldots , S_w \} $ iff
every
$P_j $ is a (not necessarily nonempty)
subset of $\{ S_1 , \ldots , S_w \} $
and the set $\{ P_1 , \ldots , P_m \} \setminus \{ \emptyset \} $ is a partition
of $\{ S_1 , \ldots , S_w \} $.

\end{define}

\begin{define}
\label{canonpart}

Let $(U,{\cal D} , G^* )$ be a standard configuration,
let $D,E\in {\cal D}$, let $C\subsetneq E$ be a
$\Lambda  ^* (E)$-block, let $Z\subseteq \Lambda  ^* (E)\cdot C$,
and
let
$S$ be a block of
$\Lambda  ^* (D)$.
The {\bf canonical  $S$-partition of $Z$}
is the
sorted
partition
$\left( S_\diamondsuit (Z), S_\| (Z), S_{\curlyvee _1 } (Z),
\ldots , S_{\curlyvee _p } (Z) \right) $
of $Z$
with the following sets.
\begin{enumerate}
\item
The set
$S_\diamondsuit (Z):=\{ A\in Z: S\diamondsuit A\} $
of blocks in $Z$ that are made by $S$.
\item
The set
$S_\| (Z):=\{ A\in Z: S\| A\} $
of blocks in $Z$ that are incomparable with $S$.
\item
The sets
$S_{\curlyvee _j } (Z):=\{ A\in Z:S\curlyvee _j A\} $
of blocks in $Z$ that are $j$-woven
with $S$.
\end{enumerate}

\end{define}

\begin{define}
\label{C-homeblock}

Let $(U,{\cal D} , G^* )$ be a standard configuration,
let
$D,E\in {\cal D}$,
let
$C$ be a block of
$\Lambda  ^* (E)$,
let $Z\subseteq \Lambda  ^* (E)\cdot C$,
and
let $S$ be a
block of
$\Lambda  ^* (D)$.
For any $S'\in \Lambda  ^* (D )\cdot S$,
we define the {\bf $Z$-home block of $S'$}, denoted
$H_{Z} (S')$, to be the set of all
$S''\in
\Lambda  ^* (D )\cdot S$
such that, for all
$\bigstar \in \{ \diamondsuit , \| ,
\curlyvee _1 , \ldots , \curlyvee _p \} $,
we have $S' _\bigstar (Z)=S'' _\bigstar (Z)$,
that is,
the canonical $S'$-partition of $Z$
equals the canonical $S''$-partition of $Z$.

\end{define}

\begin{lem}
\label{homeblocklem}

Let $(U,{\cal D} , G^* )$ be a standard configuration,
let
$D,E\in {\cal D}$, let $S$ be a
block of
$\Lambda  ^* (D)$,
and
let
$C$ be a block of
$\Lambda  ^* (E)$.
Then, for every $S'\in \Lambda  ^* (D)\cdot S$, the
$\Lambda  ^* (E)\cdot C$-home block
$H_{\Lambda  ^* (E)\cdot C} (S')$ of $S'$ is a block of
the action of
$G^* $
on $\Lambda  ^* (D)\cdot S$.

\end{lem}

{\bf Proof.}
Let $S'\in \Lambda  ^* (D)\cdot S$ and let
$\Phi \in
G^*
$
such that
$\Phi \left[ H_{\Lambda  ^* (E)\cdot C} (S')\right]
\cap H_{\Lambda  ^* (E)\cdot C} (S')\not= \emptyset $.
Because, for every $S''\in H_{\Lambda  ^* (E)\cdot C} (S')$,
we have
$H_{\Lambda  ^* (E)\cdot C} (S'')=H_{\Lambda  ^* (E)\cdot C} (S')$,
we can assume without
loss of generality that
$\Phi (S')\in H_{\Lambda  ^* (E)\cdot C} (S')$.

Let
$\bigstar \in \{ \diamondsuit , \| ,
\curlyvee _1 , \ldots , \curlyvee _p \} $
and let $C'\in {\Lambda  ^* (E)\cdot C}$.
Then
%such that
$S'\bigstar C'$
%Then
iff
$\Phi [S']\bigstar \Phi [C']$,
and, because
$\Phi (S')\in H_{\Lambda  ^* (E)\cdot C} (S')$,
%we obtain
this is the case iff
$S'\bigstar \Phi [C']$.
Hence
\begin{eqnarray*}
\lefteqn{\Phi \left[ S' _\bigstar ({\Lambda  ^* (E)\cdot C})\right]}
\\
& = &
\Phi \left[ \left\{ C'\in {\Lambda  ^* (E)\cdot C}:
S' \bigstar C'\right\} \right]
=
\left\{ \Phi \left[ C' \right] \in {\Lambda  ^* (E)\cdot C}:
S' \bigstar C'\right\}
\\
& = &
\left\{ \Phi \left[ C' \right] \in {\Lambda  ^* (E)\cdot C}:
S' \bigstar \Phi \left[ C' \right] \right\}
=
\left\{ C''\in {\Lambda  ^* (E)\cdot C}:
S' \bigstar C''\right\}
\\
& = &
S' _\bigstar \left( {\Lambda  ^* (E)\cdot C} \right)
.
\end{eqnarray*}

Consequently,
$\Phi $ maps each set in the canonical
$S'$-partition of
${\Lambda  ^* (E)\cdot C}$ to itself.
Now let
$\widetilde{S} \in H_{\Lambda  ^* (E)\cdot C} (S') $.
Then,
for every
$\bigstar \in \{ \diamondsuit , \| ,
\curlyvee _1 , \ldots , \curlyvee _p \} $,
we have
\begin{eqnarray*}
\lefteqn{\Phi \left[ \widetilde{S} \right] _\bigstar \left( {\Lambda  ^* (E)\cdot C} \right)}
\\
& = &
\left\{ \widehat{C}\in {\Lambda  ^* (E)\cdot C}:
\Phi \left[ \widetilde{S} \right] \bigstar \widehat{C}\right\}
=
\left\{ \Phi \left[ \widetilde{C}\right] \in {\Lambda  ^* (E)\cdot C}:
\Phi \left[ \widetilde{S} \right] \bigstar \Phi \left[ \widetilde{C} \right] \right\}
\\
& = &
\left\{ \Phi \left[ \widetilde{C}\right] \in {\Lambda  ^* (E)\cdot C}:
\widetilde{S} \bigstar \widetilde{C} \right\}
=
\Phi \left[ \left\{ \widetilde{C}\in {\Lambda  ^* (E)\cdot C}:
\widetilde{S} \bigstar \widetilde{C}\right\} \right]
\\
& = &
\Phi \left[ \widetilde{S} _\bigstar \left( {\Lambda  ^* (E)\cdot C} \right) \right]
=
\Phi \left[ S' _\bigstar \left( {\Lambda  ^* (E)\cdot C}\right) \right]
=
S' _\bigstar \left( {\Lambda  ^* (E)\cdot C} \right)
,
\end{eqnarray*}
that is,
$\Phi \left[ \widetilde{S} \right] \in H_{\Lambda  ^* (E)\cdot C} (S') $.
Hence
$\Phi \left[ H_{\Lambda  ^* (E)\cdot C} (S') \right]
= H_{\Lambda  ^* (E)\cdot C} (S') $
and $H_{\Lambda  ^* (E)\cdot C} (S') $ is a block of
the action of
$G^* $
on $\Lambda  ^* (D)\cdot S$.
\qed

\begin{lem}
\label{controllemma}

Let $(U,{\cal D} , G^* )$ be a standard configuration,
let
$D,E\in {\cal D}$, let $(S,B)$ be a
primitive pair for
$\Lambda  ^* (D)$,
and
let
$C$ be a block of
$\Lambda  ^* (E)$
such that
the $C$-canonical partition of
$B[\Lambda  ^* (D)\cdot S]$
has more than one nonempty set.
Then, for any $\Phi , \Psi \in G^* $,
we have that
$\Phi \hochdamit ^{\Lambda  ^* (E)\cdot C}
=\Psi \hochdamit ^{\Lambda  ^* (E)\cdot C} $
and
$\Phi \hochdamit ^{\Lambda  ^* (D)\cdot B}
=\Psi \hochdamit ^{\Lambda  ^* (D)\cdot B} $
implies that
$\Phi \hochdamit ^{\Lambda  ^* (D)\cdot S}
=\Psi \hochdamit ^{\Lambda  ^* (D)\cdot S} $.

\end{lem}

{\bf Proof.}
For any $S'\in B[\Lambda  ^* (D)\cdot S]$,
we define $K_C (S'):=H_{\Lambda  ^* (E)\cdot C} (S')
\cap B[\Lambda  ^* (D)\cdot S]$.
By Lemma \ref{homeblocklem},
$H_{\Lambda  ^* (E)\cdot C} (S')$
is a block of the action of $G^* $
on $\Lambda  ^* (D)\cdot S$.
Because the intersection of two blocks is again a block, we
have that
$K_C (S')$ is a block of the
action of $G^* $
on
$B[\Lambda  ^* (D)\cdot S]$, too.

Let $S'\in B[\Lambda  ^* (D)\cdot S]$ and
suppose, for a contradiction, that
$K_C (S')$
has more than one element.
Because
$K_C (S')$ is a block and $(S,B)$ is a primitive pair, this means that
$K_C (S')=B[\Lambda  ^* (D)\cdot S]$.
That is, for any two
$\widehat{S} , \widetilde{S} \in B[\Lambda  ^* (D)\cdot S]$
the
canonical
$\widehat{S} $-partition of
$\Lambda  ^* (E)\cdot C$ equals the
canonical
$\widetilde{S} $-partition of
$\Lambda  ^* (E)\cdot C$.
However, then
the $C$-canonical partition of
$B[\Lambda  ^* (D)\cdot S]$
has exactly one nonempty set, a contradiction.
Thus,
for every $S'\in B[\Lambda  ^* (D)\cdot S]$,
we have that $K_C (S')$
has exactly one element.

Now let
$\Phi , \Psi \in G^* $,
such that
$\Phi \hochdamit ^{\Lambda  ^* (E)\cdot C}
=\Psi \hochdamit ^{\Lambda  ^* (E)\cdot C} $
and
$\Phi \hochdamit ^{\Lambda  ^* (D)\cdot B}
=\Psi \hochdamit ^{\Lambda  ^* (D)\cdot B} $.
Because $\Phi \hochdamit ^{\Lambda  ^* (D)\cdot B}
=\Psi \hochdamit ^{\Lambda  ^* (D)\cdot B} $,
we obtain
$\Psi ^{-1} \Phi [B]=B$.
Because
$\Phi \hochdamit ^{\Lambda  ^* (E)\cdot C}
=\Psi \hochdamit ^{\Lambda  ^* (E)\cdot C} $,
we obtain
that
$
\left(
\Psi ^{-1}
\Phi
\right)
\hochdamit ^{\Lambda  ^* (E)\cdot C}
$
is the identity on $\Lambda  ^* (E)\cdot C$.
Hence,
for any $S'\in B[\Lambda  ^* (D)\cdot S]$, we have
$K_C \left( \Psi ^{-1} \Phi [S']\right ) =
K_C (S')$, which implies
$\Psi ^{-1} \Phi [S'] = S'$ and then $\Phi [S'] = \Psi [S']$.

We have thus proved that the hypotheses imply
$\Phi \hochdamit ^{\Lambda  ^* (D)\cdot S} _B
=\Psi \hochdamit ^{\Lambda  ^* (D)\cdot S} _B $.
Let $B'\in \Lambda  ^* (D)\cdot B$.
Because  $\Lambda  ^* (D)$ is transitive, there is
a $\Delta \in G^* $ such that
$\Delta [B]=B'$,
and hence
$B'$ fulfills the hypotheses with $C':=\Delta [C]$.
Via the same argument, we obtain that
$\Phi \hochdamit ^{\Lambda  ^* (D)\cdot S} _{B'}
=\Psi \hochdamit ^{\Lambda  ^* (D)\cdot S} _{B'} $.
Because
$B'\in \Lambda  ^* (D)\cdot B$ was arbitrary, we conclude
$\Phi \hochdamit ^{\Lambda  ^* (D)\cdot S}
=\Psi \hochdamit ^{\Lambda  ^* (D)\cdot S} $.
\qed

\vspace{.1in}

In light of Lemma \ref{controllemma},
we define the following.

\begin{define}
\label{controldef}

Let $(U,{\cal D} , G^* )$ be a standard configuration,
let
$D,E\in {\cal D}$, let $(S,B)$ be a
primitive pair for
$\Lambda  ^* (D)$,
and
let
$C$ be a block of
$\Lambda  ^* (E)$.
Then we define
$C \blacktriangleright ^B S $ and say that
$C $ (together with $B$)
{\bf controls}
$S $ iff
the $C$-canonical partition of
$B[\Lambda  ^* (D)\cdot S]$
has more than one nonempty set.

\end{define}

It turns out that
control will be exerted from
a larger block, too, or
control will be mutual,
as the
following lemma shows.

\begin{lem}
\label{mutualoruppercontrol}

Let $(U,{\cal D} , G^* )$ be a standard configuration,
let
$D,E\in {\cal D}$, let $(S,B)$ be a
primitive pair for
$\Lambda  ^* (D)$,
let
$C$ be a block of
$\Lambda  ^* (E)$
such that
$C \blacktriangleright ^B S $,
and
let
$T$ be a block of
$\Lambda  ^* (E)$
such that
$(C,T)$ is a primitive pair.
Then
$T \blacktriangleright ^B S $
or
$S \blacktriangleright ^T C $.

\end{lem}

{\bf Proof.}
Consider the case that
$S\not\blacktriangleright ^T C $.
Then the $S$-canonical partition of
$T[\Lambda  ^* (E)\cdot C]$
has exactly one nonempty set.
Let
$\bigstar \in \{
\diamondsuit , \| , \curlyvee _1 ,
\ldots , \curlyvee _p
\} $
such that, for all
$C'\in T[\Lambda  ^* (E)\cdot C]$,
we have $S\bigstar C'$.
Because
$C \blacktriangleright ^B S $,
there are $\widehat{S} \in
B[\Lambda  ^* (D)\cdot S ]$
and
$\blacktriangledown  \in \{
\diamondsuit , \| , \curlyvee _1 ,
\ldots , \curlyvee _p
\} \setminus \{ \bigstar \} $
such that
$\widehat{S} \blacktriangledown C$.
We conclude that
$(S , T)\not\bowtie \left( \widehat{S} , T\right) $,
which means
$T \blacktriangleright ^B S $.
\qed

\vspace{.1in}

Lemma \ref{mutualoruppercontrol} shows that control can be mutual,
and this will be an important notion in the future.
Hence, we define the following.

\begin{define}

Let $(U,{\cal D} , G^* )$ be a standard configuration,
let
$D,E\in {\cal D}$, let $(S,B)$ be a
primitive pair for
$\Lambda  ^* (D)$,
and
let
$(C,T)$ be a primitive pair for
$\Lambda  ^* (E)$.
We
say that
$(S,T)$ and $(C,T)$ are {\bf mutually controlling},
and we write
$(S,T)\directbothways (C,T)$,
iff
$C \blacktriangleright ^B S $
and $S \blacktriangleright ^T C $.

\end{define}

As Definition \ref{primnestcoll} indicates,
we will usually work with a
fixed primitive nesting collection ${\cal B}$.
When we have
a fixed primitive nesting collection ${\cal B}$,
the natural
control relations
$C\blacktriangleright ^B S$
to investigate
are those with
$B=S^\dag $.
Moreover, although
a specific block $C\in {\cal B}$ may not control
another $S\in {\cal B}$ in the exact sense of
Definition \ref{controldef},
for the
control of the group action, it suffices that
{\em some} block in
$\Lambda  ^* (E)\cdot C$ does.
Therefore, we define the following.

\begin{define}
\label{controlbyreldef}

Let $(U,{\cal D} , G^* ;{\cal B})$ be a
nested standard configuration
and let $C,S\in G^* \cdot {\cal B}$.
We will write
$C \blacktriangleright S $
instead of
$C\blacktriangleright ^{S^\dag } S $,
and we will write
$C\directbothways S$,
iff
$C \blacktriangleright S $
and $S \blacktriangleright C $.
Moreover, we will write
$C\vartriangleright S$ iff there is a conjugate
$C'$ of $C$ such that $C'\blacktriangleright S$.
We define
$C\bothways S$
iff
there is a conjugate $C'$ of $C$ such that
$C'\directbothways S$.

\end{define}

For blocks at depth zero, a controlling block always exists.
We also provide an insight into the mutual control of
blocks at depth zero that
will be needed in the proof of Proposition
\ref{breadth4Lemma}.

\begin{lem}
\label{everythingiscontrolled}

Let $(U,{\cal D} , G^* ;{\cal B} )$ be a
nested standard configuration.
Then, for every $B_0 ^D \in {\cal B}$, there is
a $B_0 ^E \in {\cal B}$ such that
$B_0 ^E \vartriangleright B_0 ^D $.

\end{lem}

{\bf Proof.}
Because $B_1 ^D $ cannot be an order-autonomous antichain,
there are $x,y\in B_1 ^D $
and a $z\not\in D$ such that
$z\sim x$ and $z\not\sim y$.
Let $E\in {\cal D}$ be the orbit that
contains $z$. Then
$B_0 ^E \vartriangleright B_0 ^D $.
\qed

\begin{lem}
\label{connectblocks}

Let $(U,{\cal D} , G^* ; {\cal B})$ be a nested standard configuration,
let
$C,S\in G^* \cdot {\cal B}$ be depth $0$ blocks
such that
$C^\dag \not\blacktriangleright S$
and
$S^\dag \not\blacktriangleright C$.
Then any two elements in
$C^\dag $ have the same number of comparable elements in
$S^\dag $, and,
any two elements in
$S^\dag $ have the same number of comparable elements in
$C^\dag $.
Moreover,
if
$C\blacktriangleright S$,
then
$C\directbothways S$
and
every element in
$C^\dag $
$\left( S^\dag {\rm , \ respectively}\right) $
is comparable to
fewer than
$\left| S^\dag \right| $
$\left( \left| C^\dag \right| {\rm , \ respectively}\right) $
elements of
$S^\dag $
$\left( C^\dag {\rm , \ respectively}\right) $.
Finally,
if
$C\blacktriangleright S$
and
$\left| S^\dag \right| $
is prime, then
$\left| S^\dag \right| $
divides
$\left| C^\dag \right| $.

\end{lem}

{\bf Proof.}
Because
$C^\dag \not\blacktriangleright S$,
for any two $S_1 , S_2 \in S^\dag [\Lambda  ^* (D_S )\cdot S]$,
we have that
$\left( C^\dag , S_1 \right) $
and $\left( C^\dag , S_2 \right) $ are similar
weavings.
Because
the depth of
$C$ and of
$S$ is $0$,
the only possible weavings for
a $C'\in C^\dag [\Lambda  ^* (D_C )\cdot C]$
with an $S'\in S^\dag [\Lambda  ^* (D_S )\cdot S]$
are the trivial weavings
$C'\lozenge S'$ or $C'\| S'$.
Therefore,
the fact that
$\left( C^\dag , S_1 \right) $ and
$\left( C^\dag , S_2 \right) $ are similar
weavings is equivalent to the fact that
the unique element of $S_1 $ is comparable to
the same number of elements in $C^\dag $ as
the unique element of $S_2 $.
Hence
any two elements in
$S^\dag $ have the same number of comparable elements in
$C^\dag $.
Symmetrically, we obtain that
any two elements in
$C^\dag $ have the same number of comparable elements in
$S^\dag $.

Let
$C\blacktriangleright S$.
Then, by Lemma \ref{mutualoruppercontrol}, because
$C^\dag \not\blacktriangleright S$,
we obtain
$S\blacktriangleright C$.
Now let
$N$ be the total number of
comparabilities between
$C^\dag $ and $S^\dag $.
Let $A_{C^\dag S^\dag } $ be the number of elements of
$S^\dag $ that are comparable to an element of $C^\dag $
(recall that this number does not depend on the
element of $C^\dag $)
and
let $A_{S^\dag C^\dag } $ be the number of elements of
$C^\dag $ that are comparable to an element of $S^\dag $.
Then
$
\left| C^\dag \right|
A_{C^\dag S^\dag }
=N=
\left| S^\dag \right|
A_{S^\dag C^\dag }
$.
Because $C\blacktriangleright S$, we have
$N<
\left| C^\dag \right|
\left| S^\dag \right|
$, which, in particular, means that
every element in
$C^\dag $ is comparable to
fewer than
$\left| S^\dag \right| $
elements of
$S^\dag $.

Now let
$\left| S^\dag \right| $
be prime.
Because
$\left| S^\dag \right| $
divides $N$,
it
divides one of
$\left| C^\dag \right| $
or
$A_{C^\dag S^\dag } $.
If
$\left| S^\dag \right| $
were to divide
$A_{C^\dag S^\dag } $,
then
$
N=
\left| C^\dag \right|
A_{C^\dag S^\dag }
\geq
\left| C^\dag \right|
\left| S^\dag \right|
$, which was excluded.
Hence,
if
$\left| S^\dag \right| $
is prime, then
$\left| S^\dag \right| $
divides
$\left| C^\dag \right| $.
\qed

\vspace{.1in}

We must end this section on a cautionary note.
Although the control relation
will be very useful indeed,
the example below shows that
controlling blocks need not exist.

\begin{figure}

\centerline{
%\input{suspended_factorial.pic}
% This is a LaTeX picture output by TeXCAD.
% File name: [suspended_factorial.pic].
% Version of TeXCAD: 4.51
% Reference / build: 27-Nov-2018 (rev. a75)
% For new versions, check: http://texcad.sf.net/
% Options on the following lines.
%\grade{\on}
%\emlines{\off}
%\epic{\off}
%\beziermacro{\on}
%\reduce{\on}
%\snapping{\on}
%\pvinsert{% Your \input, \def, etc. here}
%\quality{8.000}
%\graddiff{0.005}
%\snapasp{1}
%\zoom{10.0000}
\unitlength 1mm % = 2.845pt
\linethickness{0.4pt}
\ifx\plotpoint\undefined\newsavebox{\plotpoint}\fi % GNUPLOT compatibility
\begin{picture}(159,47)(0,0)
\put(5,15){\circle*{2}}
\put(45,15){\circle*{2}}
\put(45,5){\circle*{2}}
\put(115,35){\circle*{2}}
\put(45,35){\circle*{2}}
\put(85,15){\circle*{2}}
\put(125,15){\circle*{2}}
\put(15,15){\circle*{2}}
\put(55,15){\circle*{2}}
\put(55,5){\circle*{2}}
\put(105,35){\circle*{2}}
\put(55,35){\circle*{2}}
\put(95,15){\circle*{2}}
\put(135,15){\circle*{2}}
\put(25,15){\circle*{2}}
\put(65,15){\circle*{2}}
\put(65,5){\circle*{2}}
\put(95,35){\circle*{2}}
\put(65,35){\circle*{2}}
\put(105,15){\circle*{2}}
\put(145,15){\circle*{2}}
\put(35,15){\circle*{2}}
\put(75,15){\circle*{2}}
\put(75,5){\circle*{2}}
\put(85,35){\circle*{2}}
\put(75,35){\circle*{2}}
\put(115,15){\circle*{2}}
\put(155,15){\circle*{2}}
\put(5,25){\circle*{2}}
\put(45,25){\circle*{2}}
\put(85,25){\circle*{2}}
\put(125,25){\circle*{2}}
\put(15,25){\circle*{2}}
\put(55,25){\circle*{2}}
\put(95,25){\circle*{2}}
\put(135,25){\circle*{2}}
\put(25,25){\circle*{2}}
\put(65,25){\circle*{2}}
\put(105,25){\circle*{2}}
\put(145,25){\circle*{2}}
\put(35,25){\circle*{2}}
\put(75,25){\circle*{2}}
\put(115,25){\circle*{2}}
\put(155,25){\circle*{2}}
\put(5,15){\line(0,1){10}}
\put(45,15){\line(0,1){10}}
\put(85,15){\line(0,1){10}}
\put(125,15){\line(0,1){10}}
\put(15,15){\line(0,1){10}}
\put(55,15){\line(0,1){10}}
\put(95,15){\line(0,1){10}}
\put(135,15){\line(0,1){10}}
\put(25,15){\line(0,1){10}}
\put(65,15){\line(0,1){10}}
\put(105,15){\line(0,1){10}}
\put(145,15){\line(0,1){10}}
\put(35,15){\line(0,1){10}}
\put(75,15){\line(0,1){10}}
\put(115,15){\line(0,1){10}}
\put(155,15){\line(0,1){10}}
\put(5,25){\line(1,-1){10}}
\put(45,25){\line(1,-1){10}}
\put(85,25){\line(1,-1){10}}
\put(125,25){\line(1,-1){10}}
\put(15,25){\line(1,-1){10}}
\put(55,25){\line(1,-1){10}}
\put(95,25){\line(1,-1){10}}
\put(135,25){\line(1,-1){10}}
\put(25,25){\line(1,-1){10}}
\put(65,25){\line(1,-1){10}}
\put(105,25){\line(1,-1){10}}
\put(145,25){\line(1,-1){10}}
\put(5,15){\line(3,1){30}}
\put(45,15){\line(3,1){30}}
\put(85,15){\line(3,1){30}}
\put(125,15){\line(3,1){30}}
\put(3,13){\dashbox{1}(34,4)[cc]{}}
\put(43,13){\dashbox{1}(34,4)[cc]{}}
\put(83,13){\dashbox{1}(34,4)[cc]{}}
\put(123,13){\dashbox{1}(34,4)[cc]{}}
\put(3,23){\dashbox{1}(34,4)[]{}}
\put(43,23){\dashbox{1}(34,4)[]{}}
\put(83,23){\dashbox{1}(34,4)[]{}}
\put(123,23){\dashbox{1}(34,4)[]{}}
\put(80,15){\oval(158,6)[]}
\put(80,25){\oval(158,6)[]}
\put(5,15){\line(4,-1){40}}
\put(155,25){\line(-4,1){40}}
\put(5,25){\line(4,1){40}}
\put(15,15){\line(4,-1){40}}
\put(145,25){\line(-4,1){40}}
\put(15,25){\line(4,1){40}}
\put(25,15){\line(4,-1){40}}
\put(135,25){\line(-4,1){40}}
\put(25,25){\line(4,1){40}}
\put(35,15){\line(4,-1){40}}
\put(125,25){\line(-4,1){40}}
\put(35,25){\line(4,1){40}}
\put(45,5){\line(0,1){10}}
\put(115,35){\line(0,-1){10}}
\put(45,35){\line(0,-1){10}}
\put(55,5){\line(0,1){10}}
\put(105,35){\line(0,-1){10}}
\put(55,35){\line(0,-1){10}}
\put(65,5){\line(0,1){10}}
\put(95,35){\line(0,-1){10}}
\put(65,35){\line(0,-1){10}}
\put(75,5){\line(0,1){10}}
\put(85,35){\line(0,-1){10}}
\put(75,35){\line(0,-1){10}}
\put(85,15){\line(-4,-1){40}}
\put(75,25){\line(4,1){40}}
\put(85,25){\line(-4,1){40}}
\put(95,15){\line(-4,-1){40}}
\put(65,25){\line(4,1){40}}
\put(95,25){\line(-4,1){40}}
\put(105,15){\line(-4,-1){40}}
\put(55,25){\line(4,1){40}}
\put(105,25){\line(-4,1){40}}
\put(115,15){\line(-4,-1){40}}
\put(45,25){\line(4,1){40}}
\put(115,25){\line(-4,1){40}}
%\emline(125,15)(45,5)
\multiput(125,15)(-.2693602694,-.0336700337){297}{\line(-1,0){.2693602694}}
%\end
%\emline(35,25)(115,35)
\multiput(35,25)(.2693602694,.0336700337){297}{\line(1,0){.2693602694}}
%\end
%\emline(125,25)(45,35)
\multiput(125,25)(-.2693602694,.0336700337){297}{\line(-1,0){.2693602694}}
%\end
%\emline(135,15)(55,5)
\multiput(135,15)(-.2693602694,-.0336700337){297}{\line(-1,0){.2693602694}}
%\end
%\emline(25,25)(105,35)
\multiput(25,25)(.2693602694,.0336700337){297}{\line(1,0){.2693602694}}
%\end
%\emline(135,25)(55,35)
\multiput(135,25)(-.2693602694,.0336700337){297}{\line(-1,0){.2693602694}}
%\end
%\emline(145,15)(65,5)
\multiput(145,15)(-.2693602694,-.0336700337){297}{\line(-1,0){.2693602694}}
%\end
%\emline(15,25)(95,35)
\multiput(15,25)(.2693602694,.0336700337){297}{\line(1,0){.2693602694}}
%\end
%\emline(145,25)(65,35)
\multiput(145,25)(-.2693602694,.0336700337){297}{\line(-1,0){.2693602694}}
%\end
%\emline(155,15)(75,5)
\multiput(155,15)(-.2693602694,-.0336700337){297}{\line(-1,0){.2693602694}}
%\end
%\emline(5,25)(85,35)
\multiput(5,25)(.2693602694,.0336700337){297}{\line(1,0){.2693602694}}
%\end
%\emline(155,25)(75,35)
\multiput(155,25)(-.2693602694,.0336700337){297}{\line(-1,0){.2693602694}}
%\end
\put(60,5){\oval(34,4)[]}
\put(100,35){\oval(34,4)[]}
\put(60,35){\oval(34,4)[]}
\put(60,45){\circle*{2}}
\put(45,35){\line(3,2){15}}
\put(75,35){\line(-3,2){15}}
\put(60,45){\line(-1,-2){5}}
\put(60,45){\line(1,-2){5}}
\put(10,15){\makebox(0,0)[cc]{\footnotesize $S$}}
\put(10,25){\makebox(0,0)[cc]{\footnotesize $C$}}
\put(2,11){\makebox(0,0)[lt]{\footnotesize $B=D$}}
\put(2,29){\makebox(0,0)[lb]{\footnotesize $T=E$}}
\put(60,45){\oval(4,4)[]}
\end{picture}
}

\caption{An ordered set
$P$ with a factorial pairs $(S,B)$
and $(C,T)$ such that
neither of $C$ and $S$
is controlled by a block that is not in the union $T\cup B=E\cup D$,
and yet the weavings with these blocks are nontrivial.}
\label{suspended_factorial}

\end{figure}
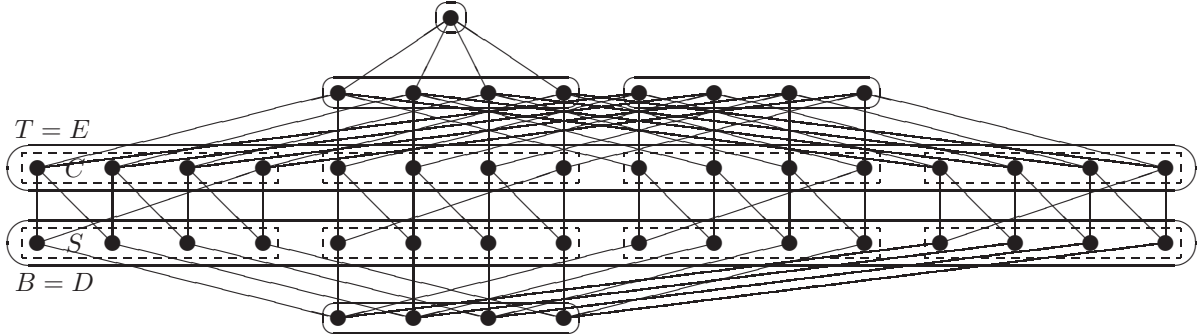

\begin{exam}
\label{suspfactex}

{\rm
In the ordered set
$P$ in Figure \ref{suspended_factorial},
equipped with its natural orbit structure,
orbits are marked with ovals and blocks under consideration are
marked with dashed boxes.
Consider the block $B=D$
and the subblock $S$.
The primitive pair $(S,B)$ is indeed a factorial pair, the action
on $B[\Lambda  ^* (D)\cdot S]$
of the subgroup of
${\rm Aut} (P)$ that fixes $B$
is the dihedral group on $4$ elements, and, for any
dashed block, any two weavings with
points in an orbit that
is not $D$ or $E$
are similar and nontrivial.
An analogous statement holds for the
block $T=E$ and the subblock $C$.
Thus, no block in $P\setminus (D\cup E)$
controls
(with the help of $B$ or $T$, respectively)
any of the blocks in
$
\Lambda  (D)\cdot S
\cup
\Lambda  (E)\cdot C
$.

This example shows how factorial actions can be
``suspended" within an ordered set
with nontrivial weavings with non-factorial
primitive pairs.
By Definition \ref{mutualmaxlockdef}
below, $(C,T)$
and
$(S,B)$
are
max-locked.
\qex
}

\end{exam}

\section{Compensation at Depth 0}
\label{depth0sec}

A
``control sequence" $C\vartriangleright S\vartriangleright Z$
need not imply that
any conjugate of $C$ controls $Z$.
Thus, in addition to the
direct inter-block combinatorial control
via the relations
$\vartriangleright $ and
$\blacktriangleright $, we also have
indirect inter-block combinatorial control
via the transitive closure of
$\vartriangleright $.
We focus our definition on depth 0.

\begin{define}
\label{depth0transitive}

Let $(U,{\cal D} , G^* ;{\cal B} )$ be a
nested standard configuration and let
${\cal B} _0 :=
\left\{ B_0 ^D : D\in {\cal D} \right\} \subset {\cal B}$.
We define $\vartriangleright _0 :=
\vartriangleright |_{{\cal B} _0 \times {\cal B} _0 } $
to be the reduction
of the
natural control relation to
${\cal B} _0 \times {\cal B} _0 $
and call it the {\bf depth 0 control relation}.
The transitive closure of
$\vartriangleright _0 $
is denoted
$\transcont _0 $.
For every $B_0 ^D \in {\cal B} _0 $,
we define
$
\left[
B_0 ^D \transcont _0
\right]
:=\{ S\in {\cal B} _0 :B_0 ^D
\transcont _0 S\} $
to be the set of all
$S\in {\cal B} _0 $
that are
{\bf (directly or indirectly) depth $0$ controlled} by
$B_0 ^D $.

\end{define}

Part \ref{level0elimlem2}
of Proposition \ref{level0elimlem}
shows that, if
$\left|
\Lambda  ^* \left( {\cal K}
\not\setcontrol B_0 ^D \right)
\right|
>0$, then the best factor
we can guarantee for the exponent is $1.54$, which
is too large.
In this section, we will
devise ways to undercut this worst-case estimate.

\begin{define}

Let $(U,{\cal D} , G^* ;{\cal B} )$ be a nested
standard configuration
and let
$B_0 ^{F } $
be a depth $0$ block.
Then we will say that
$B_0 ^{F } $
is {\bf controlled by deeper nestings}
iff
there are a
$B_\ell ^{E} \in {\cal B} $ such that $\ell \geq 1$
and a $B_0 ^{D } \in {\cal B} _0 $
such that
$B_\ell ^{E} \vartriangleright
B_0 ^{D} \transcont _0 B_0 ^{F} $.

\end{define}

So far, aside from the
mild conditions in
Definition \ref{RDEdef}, we have enjoyed complete freedom
in our choice of reverse depth enumerations.
To proceed, we will now
define opportunistic (see Definition \ref{opprtenumdef} below)
reverse depth enumerations.
Part \ref{opprtenumdef1}
of Definition \ref{opprtenumdef}
will help with a technical aspect of
the proof of
Proposition \ref{breadth4Lemma} below.
Part \ref{opprtenumdef2}
of Definition \ref{opprtenumdef}
will help in
\cite{SchAutCon2} with the generation of
``cascades" of conjugate blocks
that will help
compensate for many of the factors contributed by
blocks $S\in {\cal K} _u ^\gamma $.

\begin{define}

Let $(U,{\cal D} , G^* )$ be a standard configuration
and let $D,E\in {\cal D}$. We say that
the action of $G^* $ on
$E$ {\bf controls}
the action of $G^* $ on
$D$ iff, for all
$\Phi \in G^* $,
if $\Phi |_E ={\rm id} |_E $, then
$\Phi |_D ={\rm id} |_D $.

For an orbit $E\in {\cal D}$, we define the
{\bf controlled tail} of $E$ to be the largest set of orbits
$\{ D_0 , \ldots , D_m\} \in {\cal D}$ such that
$E=D_0$, and, for all $i\in \{ 1, \ldots , m\} $, there is a
$j\in \{ 0,\ldots , i-1\} $ such that $D_{j} $ controls $D_i $.

\end{define}

\begin{define}

Let
$\left( U,{\cal D} , G^* ; {\cal B}\right)$
be a nested
standard configuration
and let $E\in {\cal D}$. Then we say that
$E$ is {\bf cooperative} iff
$\left( B_0 ^{E } , B_1 ^{E } \right) $
is not a factorial pair of width $w\geq 5$,
and there is a $D\in {\cal D}$ such that
$|E|\leq |D|$ and
such that
the action of $G^* $ on $E$
controls the action of $G^* $ on $D$.
If
the action of $G^* $ on $D$
controls the action of $G^* $ on $E$, too,
we call $(E,D)$ a {\bf cooperating pair}.

\end{define}

Recall that,
by Lemma \ref{connectblocks}, for
any two
$S,C\in {\cal B} _0 $ that are not controlled by
deeper nestings,
we have that
$S\vartriangleright _0 C $
implies
$C\bothways S$.

\begin{define}
\label{remaindepth0}

Let
$\left( U,{\cal D} , G^* ; {\cal B}\right)$
be a nested
standard configuration.
We define ${\cal B} _0 ^r $
to be the set of all
$S\in {\cal B} _0 $ that are neither controlled by
deeper nestings,
nor in the union of the
$\transcont _0 $-closures of the depth $0$ blocks
of the orbits in the
controlled tail of a cooperating pair.

\end{define}

\begin{define}
\label{remainatdepth0}

Let
$\left( U,{\cal D} , G^* ; {\cal B}\right)$
be a nested
standard configuration.
We define ${\cal R}$ to be the
union of all
$\bothways $-components
$W$
of
$\left( {\cal B} _0 ^r , \bothways \right) $
such that, for
any two
$S,C\in W $
such that
$C\bothways S$,
the following hold.
\begin{enumerate}
\item
At least one of
$\left( S, S^\dag \right) $
and
$\left( C, C^\dag \right) $
is factorial.

\item
If $\left( C, C^\dag \right) $ is not factorial,
then
$\left| D_C \right| > \left| D_S \right| $.
\end{enumerate}
Finally, we define
${\cal G} :=
{\cal B} _0 ^r \setminus {\cal R}$.

\end{define}

\begin{define}
\label{opprtenumdef}

Let
$\left( U,{\cal D} , G^* ; {\cal B}, \overrightarrow{{\cal K}} \right)$
be a partially enumerated nested
standard configuration.
Then $\overrightarrow{\cal K}$
is called {\bf opportunistic
with control function $\gamma _{\cal K} ^* $} iff the following hold.
\begin{enumerate}
\item
\label{opprtenumdef1}
For every cooperating pair $(E,D)$,
all blocks of $E\cup D$ in ${\cal B}$ are contained in ${\cal K}$.

\item
\label{opprtenumdef2}
For all
$S\in {\cal K} _u ^\gamma \cap {\cal K} _f $,
if $S^\dag $ is the component $B_j $ of $\overrightarrow{\cal K}$, then
$S$ is the
immediately following component $B_{j+1} $.

\end{enumerate}

\end{define}

\begin{lem}
\label{opportenum}

Let
$\left( U,{\cal D} , G^* ; {\cal B}\right)$
be a nested
standard configuration.
Then there is an
opportunistic
partial reverse depth enumeration
$\overrightarrow{{\cal K}} $
of ${\cal B} \setminus {\cal G}$
with a control function $\beta _{\cal K} ^* $
such that
${\cal K} _u ^\beta \subseteq {\cal K} _f $.

\end{lem}

{\bf Proof.}
Start with $\overrightarrow{\cal K}$
being an empty listing.
Call an orbit $D$ {\bf unprocessed}
iff
${\cal K}$ contains no blocks of $D$.
Whenever multiple blocks of an orbit are added to
a current reverse depth enumeration, it shall be
understood that $S\subseteq B$ means that $B$ is added before $S$.

For the first recursive step, we will assume that,
for every processed orbit $E$ in a cooperating pair
$(E,D)$, all blocks of $E$ and all blocks of
orbits in the controlled tail of $E$
are in ${\cal K}$, and that we have a control function
$
\beta _{\cal K} ^*
$
as in
part \ref{opprtenumdef2} of Definition \ref{opprtenumdef}
such that
${\cal K} _u ^\beta \subseteq {\cal K} _f $.

Let $(D_0 , D_1 )$ be a cooperating pair
of orbits such that
$D_0 $ and $D_1 $ are unprocessed.
Let $D_2, \ldots , D_m $ be an enumeration of
the unprocessed orbits in the
controlled tail of $D_0 $ such that, for every
$j\in \{ 2, \ldots , m\} $,
there is an $i\in \{ 0,\ldots , j-1\} $
such that $D_i $ controls $D_j $.
Note that $m\geq 1$.
Extend the current
reverse depth enumeration
$\overrightarrow{\cal K}$ by adding the
blocks of the
$D_j $ such that, for $i<k$, the
blocks of $D_i $
are listed before the
blocks of $D_k $.
Because
every newly added $S$
with a successor
is added
immediately after $S^\dag $, this extension
satisfies
part \ref{opprtenumdef2} of Definition \ref{opprtenumdef}
for any extension of $\beta _{\cal K} ^* $.

We extend $
\beta _{\cal K} ^* $
as follows.
Let $C_0 \subset D_0 $ be the singleton block
of $D_0 $ in ${\cal K}$.
For every
$C\subset D_0  $, $C\in {\cal K}\setminus \{ C_0 \} $
of depth greater than $1$, we set
$
\beta _{\cal K} ^* (C)
:=
\left|
\Lambda  ^* \left( \overleftarrow{\cal K} \not\setcontrol C \right)
\right|
$.
By part \ref{level0elimlem2}
of Proposition \ref{level0elimlem}
applied to the enumeration so far,
we have
$
\prod _{C\subset D_0 , C\in {\cal K} \setminus {\cal K} _f , C\not= C_0 }
\beta _{\cal K} ^* (C)
\cdot
\left|
\Lambda  ^* \left( \overleftarrow{\cal K} \not\setcontrol C_0 \right)
\right|
=
{\Omega _{\cal K} (D_0 )
\over
\Upsilon _{\cal K} (D_0 )}
\left|
\Lambda  ^* \left( \overleftarrow{\cal K} \not\setcontrol C_0 \right)
\right|
\leq 2^{1.54 |D_0 |}
$.
We set
$\beta _{\cal K} ^* (C_0 )
:=
\max \left\{
{
\left|
\Lambda  ^* \left( \overleftarrow{\cal K} \not\setcontrol C_0 \right)
\right|
\over
2^{0.77 |D_0 |}
}
, 1\right\}
\leq
2^{0.61 |D_0 |}
$
and we obtain
$\Omega _{\cal K} (D_0 )\leq
2^{0.77|D_0 |} \Upsilon _{\cal K} (D_0 )$.

Because
$D_1 $ is controlled by $D_0 $,
for every
$S\subset D_1 $, $S\in {\cal K}$, we have that
$\left|
\Lambda  ^* \left( \overleftarrow{\cal K} \not\setcontrol S \right)
\right|
=1
$.
Let $S_0 \subset D_1 $ be the singleton block of $D_1 $ in ${\cal K}$.
For every
$S\subset D_1 $, $S\in {\cal K}\setminus \{ S_0 \} $, we set
$
\beta _{\cal K} ^* (S)
:=
\left|
\Lambda  ^* \left( \overleftarrow{\cal K} \not\setcontrol S \right)
\right|
=1
$.
Moreover,
we set
$\beta _{\cal K} ^* (S_0 )
:=
2^{0.77 |D_0 |}
$.
Hence
$\Omega _{\cal K} (D_1 )
=
2^{0.77 |D_0 |}
\leq
2^{0.77 |D_1 |}
$.
Finally,
for every
$j\in \{ 2, \ldots , m\} $, because
$D_j $ is controlled by a $D_i $ with $i<j$,
for every
$S\in {\cal K}\setminus \{ D_j \} $, we set
$
\beta _{\cal K} ^* (S)
:=
\left|
\Lambda  ^* \left( \overleftarrow{\cal K} \not\setcontrol S \right)
\right|
=1
$.
The thus obtained function is a
control function
(parts \ref{controlfunction1} and \ref{controlfunction2}
of Definition \ref{controlfunction}
are easily verified)
that is an extension
as desired
(the inequalities for $\Omega _{\cal K} (D_i )$
show that no depth $0$ blocks are added to
${\cal K} _u ^\beta $).

This process terminates with a
reverse depth enumeration
$\overrightarrow{\cal K}$
that contains all
blocks
of ${\cal B}$ whose orbits are in a cooperating pair
$(E,D)$
or in the controlled tail
of $E$
and
with a control function $\beta _{\cal K} ^* $
such that
${\cal K} _u ^\beta \subseteq {\cal K} _f $.
In particular,
part \ref{opprtenumdef1} of
Definition \ref{opprtenumdef} holds
for all extensions of this reverse depth enumeration.
In the following,
every newly added $S\in {\cal K} _u ^\beta $
%that has a successor
will be added
immediately after $S^\dag $, and hence
part \ref{opprtenumdef2} of
Definition \ref{opprtenumdef} will hold.

Let $E_1 , \ldots , E_n $ be an enumeration of all
the orbits whose depth $0$ block is in ${\cal G}$.
Add their non-depth $0$ blocks
$T$ to ${\cal K}$
such that, for $i<j$, the blocks of $E_i $ precede those of $E_j $,
and set $\beta _{\cal K} ^* (T):=
\Lambda _{\cal K} ^* \left( \overleftarrow{\cal K} \not\setcontrol T\right) $.
Because no depth $0$ blocks are added, we
maintain
${\cal K} _u ^\beta \subseteq {\cal K} _f $.

Next, let
$F_1 , \ldots , F_p $
be an
enumeration of all
the orbits whose depth $0$ block is in
${\cal R} ={\cal B} _0 ^r \setminus {\cal G}$ such that the
orbits whose depth $0$ block
is part of a factorial pair
precede all other orbits in the list.
Add their blocks
$T$ to ${\cal K}$
such that, for $i<j$, the blocks of
$F_i $ precede those of $F_j $,
and set $\beta _{\cal K} ^* (T):=
\Lambda _{\cal K} ^* \left( \overleftarrow{\cal K} \not\setcontrol T\right) $.
Note that, because the
depth $0$ blocks $S\in {\cal K} _f $
precede
the
depth $0$ blocks $S\not\in {\cal K} _f $,
we maintain ${\cal K} _u ^\beta \subseteq {\cal K} _f $.

Because all orbits with a depth $0$ block in
${\cal B} _0 ^r $ have been processed,
for all remaining orbits, the depth $0$ block is controlled by deeper nestings.

Let $H_1 , \ldots , H_q $ be an enumeration of all
the remaining orbits.
Add their non-depth $0$ blocks
$T$ to ${\cal K}$
such that, for $i<j$, the blocks of $H_i $ precede those of $H_j $,
and set $\beta _{\cal K} ^* (T):=
\Lambda _{\cal K} ^* \left( \overleftarrow{\cal K} \not\setcontrol T\right) $.
Because no depth $0$ blocks are added, we
maintain
${\cal K} _u ^\beta \subseteq {\cal K} _f $.

What remains are depth $0$ blocks that are
in
${\cal G}$ or that are controlled by deeper nestings.
Let $B_1 , \ldots , B_r $ be an enumeration of the
remaining
depth $0$ blocks
that are controlled by deeper nestings,
such that, for every $j\in \{ 1, \ldots , r\} $
there either is an $i\in \{ 1, \ldots , j-1\} $ such that
$B_i
\vartriangleright _0 B_j $,
or
there is a block $C$ of depth greater than $0$
such that
$C\vartriangleright B_j $.
Add these blocks to ${\cal K}$
in this order
and set $\beta _{\cal K} ^* (T):=
\Lambda _{\cal K} ^* \left( \overleftarrow{\cal K} \not\setcontrol T\right)
=1$.
Because the new
$\beta _{\cal K} ^* (T)$ all equal $1$, we
maintain
${\cal K} _u ^\beta \subseteq {\cal K} _f $, and
$\beta _{\cal K} ^* $ is as desired.
\qed

\begin{prop}
\label{breadth4Lemma}

Let
$\left( U,{\cal D} , G^* ; {\cal B}, \overrightarrow{{\cal K}} \right)$
be an
opportunistically partially enumerated nested
standard configuration
such that ${\cal K} $
contains
${\cal B}\setminus {\cal G} $
and
let $\beta _{\cal K} ^* $ be a control function.
Let $B_0 ^D, B_0 ^E \not\in {\cal K} $
such that
$B_0 ^E
\vartriangleright
B_0 ^D $,
neither of
$B_0 ^E $ and $B_0 ^D $ is controlled by deeper nestings,
$\left( B_0 ^{E } , B_1 ^{E } \right) $
is not a factorial pair of width $w\geq 5$,
and
$\left| E \right| \leq \left| D \right| $.
Then
there is an $F\in \{ D,E\} $
such that
there is an
opportunistic
reverse depth enumeration
$\overrightarrow{\cal L} $
with a control function
$\gamma _{\cal L} ^* $
that extends
$\beta _{\cal K} ^* $
such that
${\cal L} ={\cal K} \cup
\left[
B_0 ^{F } \transcont _0
\right] $
and ${\cal L} _u ^\beta ={\cal K}_u ^\gamma $.

\end{prop}

{\bf Proof.}
Because we are free to choose which
specific conjugate blocks to use for a primitive nesting,
we can assume without loss of generality that
$B_0 ^E \blacktriangleright B_0 ^D $.
Because
neither of
$B_0 ^E $ and $B_0 ^D $ is controlled by deeper nestings,
we have
$B_1 ^E \not\blacktriangleright B_0 ^D $
and
$B_1 ^D \not\blacktriangleright B_0 ^E $.
Hence, by
Lemma \ref{connectblocks},
we have
$B_0 ^D
\directbothways
B_0 ^E $.

For $F\in \{ D,E\} $, we can
let $m:=\left|
\left[ B_0 ^{F } \transcont _0 \right]
\setminus {\cal K} \right| $
and we can let
$
B_0 ^{F }
=: B_{1} , \ldots , B_{m} $ be an
enumeration of
$
\left[ B_0 ^{F } \transcont _0 \right]
\setminus {\cal K} $
such that
$B_{2} $ is the unique element of
$\left\{ B_0 ^D , B_0 ^E \right\}
\setminus \left\{ B_0 ^F \right\} $,
and such that,
for every $j>2 $,
there is an $i<j$ such that
$B_i \vartriangleright B_j $.
Because ${\cal B}\setminus {\cal B}_0  \subseteq {\cal K}$, clearly,
continuing
${\cal K} $
with $B_{1} , \ldots , B_{m} $
produces a reverse depth enumeration $\overrightarrow{\cal L} $.
For
every $C\in {\cal K}$, we set
$\gamma _{{\cal L}} ^* (C):=
\beta _{\cal K} ^* (C)$.
Because, for every $j\in \{ 3 , \ldots , m \} $,
we have that
${\cal K} \cup \{ B_{1} , \ldots , B_{j-1} \}
\setcontrol B_j $,
for every $j\in \{ 3 , \ldots , m \}  $, we
can set
$
\gamma _{{\cal L}} ^*
(B_j )
:=
\left|
\Lambda  ^* \left( (\overleftarrow{\cal L})
\not\setcontrol B_{j} \right)
\right|
=1$.

We will now show that, in each of a number of cases,
an
$F\in \{ D,E\} $ can be chosen such that
$
\Omega _{{\cal K} } (D)\cdot
\Omega _{{\cal K} } (E)\cdot
\left|
\Lambda  ^* \left( {\cal K}
\not\setcontrol B_0 ^F \right)
\right|
\leq
2^{0.99|D|}
2^{0.99|E|}
\cdot
\Upsilon _{{\cal K} } (D)\cdot
\Upsilon _{{\cal K} } (E)
$.
Because, with $Z$ being the unique element of
$\{ D,E\} \setminus \{ F\} $,
we have
$B_0 ^F \blacktriangleright B_0 ^Z $
and therefore $\left|
\Lambda  ^* \left( {\cal K} \cup \left\{ B_0 ^F \right\}
\not\setcontrol B_0 ^Z \right)
\right| = 1
$,
it is then possible to,
for $j\in \{ 1,2\} $,
choose
$
\gamma _{{\cal L}} ^*
(B_j )$
such that
$
\gamma _{{\cal L}} ^*
(B_j )
\prod _{T\subset D_{B_{j} } ,
T\in {\cal K} \setminus \left( {\cal K} _f \cup \left\{ B_j \right\} \right) }
\gamma _{\cal K} ^* (T)
\leq
2^{
0.99
|D_{B_{j} } | }
$.
This choice, which will not be explicitly mentioned any more,
completes the proof in each case.

For simpler notation,
for $X\in \{ D,E\} $,
we set
$
\Omega _{{\cal K} } '(X)
:=
{
\Omega _{{\cal K} } (X)
\over
\Upsilon _{{\cal K} } (X)
} $.
Because
$\left( B_0 ^{E } , B_1 ^{E } \right) $
is not a factorial pair of width $w\geq 5$,
we have
$\left|
\Lambda  ^* \left( {\cal K}
\not\setcontrol B_0 ^E \right)
\right|
\leq
2^{d_{\left| B_{1} ^E \right| /
\left| B_{0} ^E \right| } |E|}
$.
Via part \ref{level0elimlem2} of
Proposition \ref{level0elimlem},
we obtain
$\Omega _{{\cal K} } '(E)\cdot
\left|
\Lambda  ^* \left( {\cal K}
\not\setcontrol B_0 ^E \right)
\right|
\leq
2^{1.54 |E|}
$.

{\em Case 1: $\left| B_1 ^D \right| \geq 4$.}
Because
$\left| B_1 ^D \right| \geq 4$,
by
part \ref{level0elimlem1} of
Proposition \ref{level0elimlem},
we have $\Omega _{{\cal K} }'
(D)
\leq
2^{1.54 {|D|\over 4} }
$.
Hence the following holds.
\begin{eqnarray*}
\lefteqn{
\Omega _{{\cal K} } '(D)\cdot
\Omega _{{\cal K} } '(E)\cdot
\left|
\Lambda  ^* \left( {\cal K}
\not\setcontrol B_0 ^E \right)
\right|
}\\
& \leq &
2^{1.54 |E|}
2^{1.54 {|D|\over 4} }
\leq
2^{0.96 |E|+0.58|D|}
2^{0.39 |D| }
\leq
2^{0.96 |E|}
2^{0.97 |D| }
\end{eqnarray*}

{\em Case 2: $\left| B_1 ^D \right| =3$.}
Because
$\left| B_1 ^D \right| =3$,
by
part \ref{level0elimlem1} of
Proposition \ref{level0elimlem},
we have
$\Omega _{{\cal K} } '
(D)
\leq
2^{1.54 {|D|\over 3} }
$, and,
by Lemma \ref{connectblocks},
$\left| B_1 ^E \right| $ is a multiple of $3$.

{\em Case 2.1:
$|E|\leq {3\over 4} |D|$.}
In this case, we obtain
\begin{eqnarray*}
\lefteqn{
\Omega _{{\cal K} } '(D)\cdot
\Omega _{{\cal K} } '(E)\cdot
\left|
\Lambda  ^* \left( {\cal K}
\not\setcontrol B_0 ^E \right)
\right|
}
\\
& \leq &
2^{1.54 |E|}
2^{1.54 {|D|\over 3} }
\leq
2^{0.96 |E|+0.58\cdot {3\over 4} |D|}
2^{0.52 |D| }
\leq
2^{0.96 |E|}
2^{0.97 |D| }
\end{eqnarray*}

{\em Case 2.2:
$\left| B_1 ^E \right| =3$.}
When the action on the
blocks
$\Lambda  ^* (E)\cdot B_1 ^E $ is the identity, the
action on the
singletons in
$\Lambda  ^* (E)\cdot B_0 ^E $ decomposes into
${|E|\over 3} $ subsets on which we have at most $3!=6$ permutations each.
Hence, by
part \ref{level0elimlem1} of
Proposition \ref{level0elimlem},
%because $|E|\geq 3$,
we obtain $\Omega _{{\cal K} }'
(E)
\cdot
\left|
\Lambda  ^* \left( {\cal K}
\not\setcontrol B_0 ^E \right)
\right|
\leq
2^{1.54 {|E|\over 3} }
6^{|E|\over 3}
$
and
we obtain the following estimate.
\begin{eqnarray*}
\lefteqn{
\Omega _{{\cal K} } '(D)\cdot
\Omega _{{\cal K} } '(E)\cdot
\left|
\Lambda  ^* \left( {\cal K}
\not\setcontrol B_0 ^E \right)
\right|
}
\\
& \leq &
2^{1.54 {|E|\over 3} }
6^{|E|\over 3}
2^{1.54 {|D|\over 3} }
\leq
2^{0.52 |E|}
2^{0.87 |E|}
2^{0.52 |D| }
\leq
2^{0.96 |E|}
2^{0.95 |D| }
\end{eqnarray*}

{\em Case 2.3:
$|E|> {3\over 4} |D|$
and
$\left| B_1 ^E \right| >3$.}
Because
$\left| B_1 ^E \right| $ is a multiple of $3$,
we have $\left| B_1 ^E \right| \geq 6$
and hence, by
part \ref{level0elimlem1} of
Proposition \ref{level0elimlem},
$\Omega _{{\cal K} } '(E)
\leq
2^{{1.54\over 6} |E|}
$.
Because $\left| B_1 ^D \right| =3$,
$\left( B_0 ^D , B_1 ^D \right) $
clearly is not a factorial pair of width
$\geq 5$.
Hence
$\left|
\Lambda  ^* \left( {\cal D}
\not\setcontrol B_0 ^D \right)
\right|
\leq
2^{d_{\left| B_{1} ^D \right| /
\left| B_{0} ^D \right| } |D|}
$,
and,
via
part \ref{level0elimlem2} of
Proposition \ref{level0elimlem},
we obtain
$\Omega _{{\cal K} } '(D)\cdot
\left|
\Lambda  ^* \left( {\cal K}
\not\setcontrol B_0 ^D \right)
\right|
\leq
2^{1.54 |D|}
$,
and then
the following holds.
\begin{eqnarray*}
\lefteqn{
\Omega _{{\cal K} } '(D)\cdot
\Omega _{{\cal K} } '(E)\cdot
\left|
\Lambda  ^* \left( {\cal K}
\not\setcontrol B_0 ^D \right)
\right|
}\\
& \leq &
2^{1.54 |D|}
2^{1.54 {|E|\over 6} }
\leq
2^{0.99 |D| + 0.55\cdot {4\over 3} |E| }
2^{{1.54\over 6} |E|}
\leq
2^{0.99 |D|  }
2^{{4.4+1.54\over 6} |E|}
=
2^{0.99 |D| }
2^{0.99 |E|}
\end{eqnarray*}

The above concludes the proof for {\em Case 2}.

{\em Case 3: $\left| B_1 ^D \right| =2$.}
By Lemma \ref{connectblocks},
for every $z\in B_1 ^{E } $, we have that
$z$ is comparable to one element of
$B_1 ^{D } $, but not the other.
By transitivity of
$\Lambda  ^* (D )$, for every
block
$S\in \Lambda  ^* (D)\cdot B_1 ^{D }$,
there is a block $C\in \Lambda  ^* (E)\cdot B_1 ^{E }$
such that,
for any $z\in C $, we have that
$z$ is comparable to one element of
$S$, but not the other.
For every block
$C\in \Lambda  ^* (E)\cdot B_1 ^{E }$,
pick an element $e_C \in C$.

Let
$\Phi , \Psi \in G^* $ such that
$
\Phi \hochdamit ^{\Lambda  ^* (D )\cdot B_1 ^{D } }
=
\Psi \hochdamit ^{\Lambda  ^* (D )\cdot B_1 ^{D } }
$
and, for all
$C\in \Lambda  ^* (E)\cdot B_1 ^{E }$,
we have $\Phi (e_C )=\Psi (e_C )$.
Let
$\{ a,b\} \in \Lambda  ^* (D)\cdot B_1 ^{D }$
and let
$C\in \Lambda  ^* (E)\cdot B_1 ^{E }$
such that, without loss of generality,
$e_C \sim a$ and $e_C \not\sim b$.
Because
$
\Phi \hochdamit ^{\Lambda  ^* (D )\cdot B_1 ^{D } }
=
\Psi \hochdamit ^{\Lambda  ^* (D )\cdot B_1 ^{D } }
$,
we obtain
$\Phi [\{ a,b\} ]=\Psi [\{ a,b\} ]$.
Because
$e_C \sim a$,
we obtain
$\Phi (a)\sim \Phi (e_C )=\Psi (e_C ) \sim \Psi (a)$, and,
because
$e_C \not\sim b$,
we obtain
$\Phi (b)\not\sim \Phi (e_C )=\Psi (e_C )\not\sim \Psi (b)$.
Together with
$\Phi [\{ a,b\} ]=\Psi [\{ a,b\} ]$,
we obtain
$\Phi (a )=\Psi (a )$
and
$\Phi (b )=\Psi (b )$.
That is, for every
$\Phi \in G^* $,
$
\Phi \hochdamit ^{\Lambda  ^* (D )\cdot B_0 ^{D } }
$
is determined by
$
\Phi \hochdamit ^{\Lambda  ^* (D )\cdot B_1 ^{D } }
$ and
$\Phi |_{\left\{ e_C:C\in \Lambda  ^* (E)\cdot B_1 ^{E }\right\} } $.
Because
$
\Phi \hochdamit ^{\Lambda  ^* (D )\cdot B_0 ^{D } }
$
and
$
\Phi \hochdamit ^{\Lambda  ^* (E )\cdot B_1 ^{E } }
$
determine
$
\Phi \hochdamit ^{\Lambda  ^* (E )\cdot B_0 ^{E } }
$, we obtain that
$\Phi |_{D\cup E} $ is
determined by
$
\Phi \hochdamit ^{\Lambda  ^* (E )\cdot B_1 ^{E } }
$,
$
\Phi \hochdamit ^{\Lambda  ^* (D )\cdot B_1 ^{D } }
$
and
$\Phi |_{\left\{ e_C:C\in \Lambda  ^* (E)\cdot B_1 ^{E }\right\} } $.

Let $s:=\left| B_1 ^{E } \right| $.
When the action on the
$\Lambda  ^* (E)\cdot B_1 ^E $ is the identity,
for every $C\in \Lambda  ^* (E)\cdot B_1 ^E $,
there
are exactly $s$ possible choices for the image of $e_C $.
Therefore, by
part \ref{level0elimlem1} of
Proposition \ref{level0elimlem}, the following holds.
\begin{eqnarray*}
\lefteqn{
\Omega _{{\cal K} } '(D)\cdot
\Omega _{{\cal K} } '(E)\cdot
\left|
\Lambda  ^* \left( {\cal K}
\not\setcontrol B_0 ^E \right)
\right|
}\\
& \leq &
2^{1.54 \cdot {1\over 2} |D |}
\cdot
2^{1.54 \cdot {1\over s} |E |}
\cdot
s^{|E| \over s}
=
2^{0.77 |D |}
\cdot
2^{{1.54 \over s} |E |}
\cdot
2^{{\lg (s) \over s}|E|}
\end{eqnarray*}

By Lemma \ref{connectblocks},
$s$ is a multiple of $2$.

{\em Case 3.1:
$s\geq 4$.}
In this case, we have ${\lg (s)\over s} \leq {\lg (4)\over 4} ={1\over 2} $.
Hence,
\begin{eqnarray*}
\lefteqn{
\Omega _{{\cal K} } '(D)\cdot
\Omega _{{\cal K} } '(E)\cdot
\left|
\Lambda  ^* \left( {\cal K}
\not\setcontrol B_0 ^E \right)
\right|
}\\
& \leq &
2^{0.77 |D |}
\cdot
2^{{1.54 \over 4} |E |}
\cdot
2^{{1 \over 2}|E|}
\leq
2^{0.77 |D |}
\cdot
2^{0.39 |E |}
\cdot
2^{{1 \over 2}|E|}
\leq
2^{0.77 |D |}
\cdot
2^{0.89 |E |}
\end{eqnarray*}

{\em Case 3.2: $s=2$ and
$|E|\leq {1\over 2} |D|$.}
In this case, we have
\begin{eqnarray*}
\lefteqn{
\Omega _{{\cal K} } '(D)\cdot
\Omega _{{\cal K} } '(E)\cdot
\left|
\Lambda  ^* \left( {\cal K}
\not\setcontrol B_0 ^E \right)
\right|
}\\
& \leq &
2^{0.77|D| }
\cdot
2^{{1.54 \over 2}|E| }
\cdot
2^{ {\lg(2)\over 2} |E|}
\leq
2^{0.77 |D| }
\cdot
2^{0.43 |E|+0.34\cdot {1\over 2} |D|}
\cdot
2^{ {1\over 2} |E|}
\leq
2^{0.94 |D| }
\cdot
2^{0.93 |E|}
\end{eqnarray*}

{\em Case 3.3: $s=2$ and
$|E|> {1\over 2} |D|$.}
By Lemma \ref{connectblocks},
because
$\left| B_1 ^D \right| =\left| B_1 ^E \right| =2$, for
any
$C\in \Lambda  ^* (E)\cdot B_1 ^E $
and
$S\in \Lambda  ^* (D)\cdot B_1 ^D $
such that $C$ and $S$ are nontrivially woven, which we denote
$C\curlyvee S$, we have that
$C\cup S$ is
the disjoint union of two chains with $2$ elements each.

First suppose, for a contradiction,
that,
for every
$C\in \Lambda  ^* (E)\cdot B_1 ^E $,
there is exactly one
$S\in \Lambda  ^* (D)\cdot B_1 ^D $
such that
$C\curlyvee S$.
By transitivity,
every
$S\in \Lambda  ^* (D)\cdot B_1 ^D $
is nontrivially woven with the same number of
blocks
$C\in \Lambda  ^* (E)\cdot B_1 ^E $.
By assumption, for any two distinct
$S_1 , S_2 \in \Lambda  ^* (D)\cdot B_1 ^D $, we have that
$
\{ C\in \Lambda  ^* (E)\cdot B_1 ^E :C\curlyvee S_1 \}
\cap
\{ C\in \Lambda  ^* (E)\cdot B_1 ^E :C\curlyvee S_2 \}
=\emptyset $.
Therefore,
existence of an
$S\in \Lambda  ^* (D)\cdot B_1 ^D $
that is nontrivially woven
with $k\geq 2$
blocks
$C\in \Lambda  ^* (E)\cdot B_1 ^E $,
would lead to
$\left| \Lambda  ^* (E)\cdot B_1 ^E \right|
=k\cdot
\left| \Lambda  ^* (D)\cdot B_1 ^D \right|
$
and then
$|E|=2\cdot k\cdot {\left| D \right| \over \left| B_1 ^D \right| }
=k|D|>|D|$,
a contradiction to
$|E|\leq |D|$.
Hence, in this case,
the $\curlyvee $-relation
is a one-to-one matching between the blocks
in $\Lambda  ^* (E)\cdot B_1 ^E $
and the blocks in $\Lambda  ^* (D)\cdot B_1 ^D $.
Because every $\Phi \in G^* $ preserves $\curlyvee $,
the action of $G^* $
on
$\Lambda  ^* (E)\cdot B_1 ^E $
determines
the action of $G^* $
on
$\Lambda  ^* (D)\cdot B_1 ^D $.
Consequently,
the action of $G^* $
on
$E $
determines
the action of $G^* $
on
$D $.
Similarly,
the action of $G^* $
on
$D $
determines
the action of $G^* $
on
$E $, that is, $(E,D)$ is a cooperating pair.
Because $\overrightarrow{\cal K}$ was opportunistic,
all ${\cal B} $-blocks of
$D$ and $E$ are already in ${\cal K}$,
a contradiction.

We are left with the case that
$s=2$,
$|E|> {1\over 2} |D|$, and,
for every
$C\in \Lambda  ^* (E)\cdot B_1 ^E $,
there are at least two distinct
$S_1 , S_2 \in \Lambda  ^* (D)\cdot B_1 ^D $
such that $C\curlyvee S_j $.
Because $|E|>{1\over 2} |D|$,
by the pigeonhole principle
there is an
$S\in \Lambda  ^* (D)\cdot B_1 ^D $
such that
there are at least two distinct
$C_1 , C_2 \in \Lambda  ^* (E)\cdot B_1 ^E $
such that $C_i \curlyvee S $.
Because each $C_i \cup S$ is the union of two
disjoint chains with 2 elements each,
we have that
$C_1 \cup S\cup C_2 $ is the
disjoint union of two ordered sets
$\{ a<b>c\} $ or $\{ a>b<c\} $.
Hence, the action of $G^* $ on
$C_1 $, together with the
action on $\Lambda  ^* (E)\cdot B_1 ^E $, determines
the action of $G^* $ on
$C_2 $.

Because the action of $G^* $ is transitive on the orbits,
for every
$C\in \Lambda  ^* (E)\cdot B_1 ^E $,
there is at least one
$C'\in \Lambda  ^* (E)\cdot B_1 ^E \setminus \{ C\} $
such that the
action of $G^* $ on $C$, together with the
action on $\Lambda  ^* (E)\cdot B_1 ^E $,
determines the
action of $G^* $ on $C'$.
Because $|C|=2$, when the action on
$\Lambda  ^* (E)\cdot B_1 ^E $
is the identity, $G^* $ can induce at most
2 permutations on $C\cup C'$.

Hence $E$ decomposes into at most ${|E|\over 4} $ subsets
on which, when the action on
$\Lambda  ^* (E)\cdot B_1 ^E $
is the identity, we have 2 permutations each.
Therefore, once again by
part \ref{level0elimlem1} of
Proposition \ref{level0elimlem},
we obtain the following.
\begin{eqnarray*}
\lefteqn{
\Omega _{{\cal K} } '(D)\cdot
\Omega _{{\cal K} } '(E)\cdot
\left|
\Lambda  ^* \left( {\cal K}
\not\setcontrol B_0 ^E \right)
\right|
}\\
& \leq &
2^{1.54 {|D|\over 2} }
2^{1.54 {|E|\over 2} }
2^{|E|\over 4}
\leq
2^{0.77 |D|}
2^{0.77 |E|}
2^{0.13 |E|+ 0.12 |D| }
\leq
2^{0.89 |D| }
2^{0.90 |E|}
\end{eqnarray*}

\vspace*{-.25in}
~\qed

\begin{theorem}
\label{upperboundwithfactorials}

Let
$\left( U,{\cal D} , G^* ; {\cal B}\right)$
be a nested
standard configuration.
Then there is an
opportunistically
enumerated nested
standard configuration
$\left( U,{\cal D} , G^* ; \overrightarrow{{\cal B}} , {\cal F} \right)$
with a control function $\gamma _{\cal B} ^* $
such that
$\displaystyle{
\left| G^* \right|
\leq
2^{0.99 |U|}
\prod _{D\in {\cal B} }
\Upsilon _{{\cal B} } (D)
, } $
and
${\cal B} _u ^\gamma \subseteq {\cal B} _f $.

\end{theorem}

{\bf Proof.}
Starting with the
opportunistically partially enumerated nested
standard configuration
$\left( U,{\cal D} , G^* ; {\cal B}, \overrightarrow{{\cal K}} \right)$
from Lemma \ref{opportenum},
repeatedly apply Proposition \ref{breadth4Lemma}
until the
enumeration has reached all of ${\cal B}$.
The inequality follows from
Proposition \ref{reindexed}.
\qed

\vspace{.1in}

The upshot of
Theorem \ref{upperboundwithfactorials} is that,
if we can control factorial pairs, be it
directly as in
Definitions \ref{controldef} and
\ref{controlbyreldef}, or transitively, similar to
Definition \ref{depth0transitive}, but at arbitrary depths,
or, if we can compensate for them
(possibly using blocks $S\in {\cal R}$ such that
$\left( S, S^\dag \right) $ is not factorial
but necessarily $\directbothways $-linked with a factorial pair),
then we can obtain a bound of
$|{\rm Aut}  (U)|
\leq
2^{0.99 |U|} $.
Therefore, for the
future,
the focus lies squarely
on factorial pairs.

\section{Max-Locked Factorial Pairs}
\label{factpairinteract}

The max-lock relation $\blacksquare $, see
Definition \ref{mutualmaxlockdef} below,
is the simplest notion of mutual control between
primitive pairs.
(Remark \ref{outlook} provides an outlook for the other
notions.)
Sets of max-locked factorial pairs will be
our main focus for the remainder of this presentation.
They are called
``max-locked" because of the ``transmission
property" in Lemma \ref{thetatransmission}
below, which
is
crucial for
our argument,
and which will be expanded further in
Lemma \ref{mutualmaxlockcontrol} below.

\begin{define}
\label{mutualmaxlockdef}

Let $(U,{\cal D} , G^* )$ be a standard configuration,
let $D,E\in {\cal D}$,
let $w\geq 5$,
let
$(S ,B )$ be a
primitive
pair of width $w$
for
$\Lambda  ^* (D)$,
and let $(C,T)$ be a primitive pair of width $w$
for
$\Lambda  ^* (E)$.
We say that
$(C,T)$
and
$(S,B)$
are {\bf max-locked}
and write
$(S,B) \blacksquare (C,T)$ iff
$T[\Lambda  ^* (E)\cdot C]=\{ C_1 , \ldots , C_w \} $,
$B[\Lambda  ^* (D)\cdot S]=\{ S_1 , \ldots , S_w \} $
and
there are
$\bigstar ,\blacktriangledown \in \{
\diamondsuit , \| , \curlyvee _1 ,
\ldots , \curlyvee _p
\} $
such that, for every $j\in \{ 1,\ldots , w\} $, we have
$S_j \bigstar C_j $, and,
for every $i\in \{ 1,\ldots , w\} \setminus \{ j\} $,
we have
$S_i \blacktriangledown C_j $.
Moreover, we define
$\Theta _{B,T} :
B[\Lambda  ^* (D)\cdot S]
\to
T[\Lambda  ^* (E)\cdot C]$
by
$\Theta _{B,T} (S_j ):=C_j $
and
we define
$\Theta _{T,B} :
T[\Lambda  ^* (E)\cdot C]
\to
B[\Lambda  ^* (D)\cdot S]
$
by $\Theta _{T,B} (C_j ):= S_j $.

\end{define}

\begin{lem}
\label{thetatransmission}

Let $(U,{\cal D} , G^* )$ be a standard configuration,
let $w\geq 5$,
let $D, E\in {\cal D}$,
let $(C,T)$ be a primitive pair
of width $w$
for
$\Lambda  ^* (E)$,
and
let
$(S ,B )$ be a
primitive pair
of width $w$
for
$\Lambda  ^* (D)$
such that
$(S,B) \blacksquare (C,T)$.
Then,
for all $\Phi \in G^* $, we have
$\Phi
\hochdamit ^{\Lambda  ^* (E)\cdot C}
\circ \Theta _{B,T}
=
\left.
\Theta _{\Phi [B],\Phi [T]} \circ \Phi
\hochdamit ^{\Lambda  ^* (D)\cdot \Phi [S]}
\right| _{B[\Lambda  ^* (D)\cdot S]}
$.

\end{lem}

{\bf Proof.}
First note that,
because $\Phi \in G^* $ is an automorphism, for any
$
\blacktriangle
\in \{
\bigstar ,\blacktriangledown
\} $ and any $i,j\in \{ 1, \ldots , w\} $,
we have
$\Phi [S_i ]
\blacktriangle
\Phi [C_j ]$ iff
$S_i
\blacktriangle
C_j $.
Hence the ``image block sets"
$\Phi [T][\Lambda  ^* (E)\cdot \Phi [C]]$
and $\Phi [B][\Lambda  ^* (D)\cdot \Phi [S]]$
are max-locked
and, for all $i\in \{ 1, \ldots , w\} $, we have
$\Theta _{\Phi [B],\Phi [T]}
(\Phi [S_i ])=\Phi [C_i ]=
\Phi [\Theta _{B,T} (S_i )]$, which
establishes the equation.
\qed

\vspace{.1in}

To explore how the ``transmission property"
from Lemma \ref{thetatransmission} can be applied to
factorial pairs that are not max-locked,
but connected
%to each other
via $\blacksquare $,
we introduce the factorial max-lock graph
and the idea of synchronicity.
%in Definition \ref{synchcompdef} below.

\begin{define}
\label{synchcompdef}

Let $(U,{\cal D} , G^* ;{\cal B} , {\cal F})$ be a
nested standard configuration.
The graph
$({\cal F},\blacksquare )$
is called the
{\bf ($
{\cal B}$-)factorial max-lock
graph}.
A
$\blacksquare $-connected subset ${\cal S}\subseteq {\cal F}$ is called
{\bf synchronized} iff
we have that, for all
max-lock cycles
$
(S_1 , B_1 )\blacksquare
(S_2 , B_2 )\blacksquare
\cdots \blacksquare
(S_n , B_n )\blacksquare
(S_1 , B_1 )$
with repeated vertices allowed,
we have that
$
\Theta _{B_{n} , B_1 }
\circ
\Theta _{B_{n-1} , B_{n} }
\circ
\cdots
\circ
\Theta _{B_{1} , B_2 }
={\rm id} \hochdamit ^ {\Lambda  ^* (D_{S_1 } )\cdot S_1 } _{B_1}
$.
A $\blacksquare $-connected subset ${\cal U}\subseteq {\cal F}$
that is not synchronized is called
{\bf unsynchronized}.

\end{define}

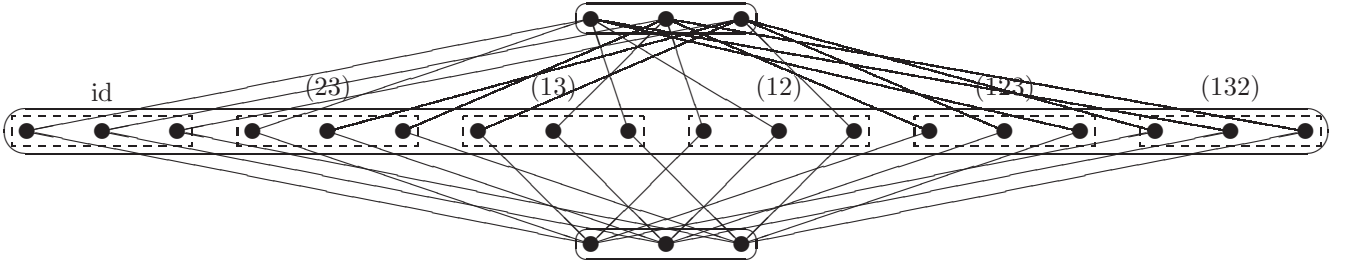
\begin{figure}

\centerline{
%\input{unsynch_factorial.pic}
% This is a LaTeX picture output by TeXCAD.
% File name: [unsynch_factorial.pic].
% Version of TeXCAD: 4.51
% Reference / build: 27-Nov-2018 (rev. a75)
% For new versions, check: http://texcad.sf.net/
% Options on the following lines.
%\grade{\on}
%\emlines{\off}
%\epic{\off}
%\beziermacro{\on}
%\reduce{\on}
%\snapping{\off}
%\pvinsert{% Your \input, \def, etc. here}
%\quality{8.000}
%\graddiff{0.005}
%\snapasp{1}
%\zoom{10.0000}
\unitlength 1mm % = 2.845pt
\linethickness{0.4pt}
\ifx\plotpoint\undefined\newsavebox{\plotpoint}\fi % GNUPLOT compatibility
\begin{picture}(178,37)(0,0)
\put(5,20){\circle*{2}}
\put(35,20){\circle*{2}}
\put(80,5){\circle*{2}}
\put(80,35){\circle*{2}}
\put(65,20){\circle*{2}}
\put(95,20){\circle*{2}}
\put(125,20){\circle*{2}}
\put(155,20){\circle*{2}}
\put(15,20){\circle*{2}}
\put(45,20){\circle*{2}}
\put(90,5){\circle*{2}}
\put(90,35){\circle*{2}}
\put(75,20){\circle*{2}}
\put(105,20){\circle*{2}}
\put(135,20){\circle*{2}}
\put(165,20){\circle*{2}}
\put(25,20){\circle*{2}}
\put(55,20){\circle*{2}}
\put(100,5){\circle*{2}}
\put(100,35){\circle*{2}}
\put(85,20){\circle*{2}}
\put(115,20){\circle*{2}}
\put(145,20){\circle*{2}}
\put(175,20){\circle*{2}}
\put(3,18){\dashbox{1}(24,4)[cc]{}}
\put(33,18){\dashbox{1}(24,4)[cc]{}}
\put(63,18){\dashbox{1}(24,4)[cc]{}}
\put(93,18){\dashbox{1}(24,4)[cc]{}}
\put(123,18){\dashbox{1}(24,4)[cc]{}}
\put(153,18){\dashbox{1}(24,4)[cc]{}}
\put(90,20){\oval(176,6)[]}
\put(90,5){\oval(24,4)[]}
\put(90,35){\oval(24,4)[]}
\put(80,5){\line(-5,1){75}}
\put(80,35){\line(-5,-1){75}}
\put(100,5){\line(5,1){75}}
\put(90,5){\line(-5,1){75}}
\put(90,35){\line(-5,-1){75}}
\put(90,5){\line(5,1){75}}
\put(100,5){\line(-5,1){75}}
\put(100,35){\line(-5,-1){75}}
\put(80,5){\line(5,1){75}}
\put(35,20){\line(3,-1){45}}
\put(35,20){\line(3,1){45}}
\put(145,20){\line(-3,-1){45}}
\put(45,20){\line(3,-1){45}}
\put(135,20){\line(-3,-1){45}}
\put(55,20){\line(3,-1){45}}
\put(125,20){\line(-3,-1){45}}
\put(65,20){\line(1,-1){15}}
\put(115,20){\line(-1,-1){15}}
\put(115,20){\line(-1,1){15}}
\put(75,20){\line(1,-1){15}}
\put(75,20){\line(1,1){15}}
\put(105,20){\line(-1,-1){15}}
\put(85,20){\line(1,-1){15}}
\put(95,20){\line(-1,-1){15}}
%\emline(45,20)(100,35)
\multiput(45,20)(.1235955056,.0337078652){445}{\line(1,0){.1235955056}}
%\end
%\emline(90,35)(55,20)
\multiput(90,35)(-.0786516854,-.0337078652){445}{\line(-1,0){.0786516854}}
%\end
%\emline(65,20)(100,35)
\multiput(65,20)(.0786516854,.0337078652){445}{\line(1,0){.0786516854}}
%\end
\put(85,20){\line(-1,3){5}}
\put(95,20){\line(-1,3){5}}
\put(105,20){\line(-5,3){25}}
%\emline(125,20)(90,35)
\multiput(125,20)(-.0786516854,.0337078652){445}{\line(-1,0){.0786516854}}
%\end
%\emline(100,35)(135,20)
\multiput(100,35)(.0786516854,-.0337078652){445}{\line(1,0){.0786516854}}
%\end
%\emline(145,20)(80,35)
\multiput(145,20)(-.1460674157,.0337078652){445}{\line(-1,0){.1460674157}}
%\end
%\emline(155,20)(100,35)
\multiput(155,20)(-.1235955056,.0337078652){445}{\line(-1,0){.1235955056}}
%\end
%\emline(165,20)(80,35)
\multiput(165,20)(-.191011236,.0337078652){445}{\line(-1,0){.191011236}}
%\end
%\emline(175,20)(90,35)
\multiput(175,20)(-.191011236,.0337078652){445}{\line(-1,0){.191011236}}
%\end
\put(15,24){\makebox(0,0)[cb]{\footnotesize ${\rm id} $}}
\put(45,24){\makebox(0,0)[cb]{\footnotesize $(23)$}}
\put(75,24){\makebox(0,0)[cb]{\footnotesize $(13) $}}
\put(105,24){\makebox(0,0)[cb]{\footnotesize $(12)$}}
\put(135,24){\makebox(0,0)[cb]{\footnotesize $(123) $}}
\put(165,24){\makebox(0,0)[cb]{\footnotesize $(132) $}}
\end{picture}
}

\caption{
An unsynchronized
subset of a graph $({\cal F} , \blacksquare )$
with, for a smaller image, $w=3$.
The 3-point dashed blocks are created via comparabilities that are
not shown.
Every combination of permutations of the top and bottom
orbit is induced by an automorphism:
If $\sigma $ is a permutation of the
top orbit and the identity acts on the bottom orbit,
assume that every block's vertices are labeled $1,2,3$ from left to right.
Then every 3-point block
$\mu $
in the middle orbit is mapped to $\sigma \mu $.
}
\label{unsynch_factorial}

\end{figure}

Figure \ref{unsynch_factorial}
shows
an example of an unsynchronized
subset of ${\cal F}$ with, by slight abuse of
notation, $w=3$.
We will frequently refer to Figure \ref{unsynch_factorial}
to illustrate the necessity of a condition or the
effect of an action.
It is also worth mentioning that the
unsynchronized factorials in
Figure \ref{unsynch_factorial}
led to the group representation in \cite{SchPermRep}

\begin{define}

Let $(U,{\cal D} , G^* ;{\cal B} , {\cal F})$ be a
nested standard configuration.
A subset ${\cal P}\subseteq {\cal F}$ is called
{\bf saturated} iff it is
a component of
$G^* \cdot {\cal P}
$
with adjacency relation $\blacksquare $.
Two saturated subsets ${\cal L} ,{\cal P}$ of
$({\cal F},\blacksquare )$ will be called
{\bf conjugate} and we will write
${\cal L}\asymp {\cal P}$,
iff
${\cal L}\subseteq G^* \cdot {\cal P}$.

\end{define}

\begin{lem}
\label{longthetawelldef}

Let $(U,{\cal D} , G^* ;{\cal B} ,{\cal F})$ be a
nested standard configuration.
Let the
saturated
subset ${\cal S}\subseteq {\cal F}$ be
synchronized and let
$(C,T), (S, B )\in {\cal S}$.
Then, for any two paths
$(C,T)=(S_0 , B_0 )\blacksquare
(S_1 , B_1 )\blacksquare
\cdots \blacksquare
(S_m , B_m )=
(S , B )$
and
$(C,T)=(S_0 ', B_0 ')\blacksquare
(S_1 ', B_1 ')\blacksquare
\cdots \blacksquare
(S_n ', B_n ')=
(S , B )$, we have
$
\Theta _{B_{m-1} , B_m }
\circ
\cdots
\circ
\Theta _{B_{0} , B_1 }
=
\Theta _{B_{n-1} ', B_n '}
\circ
\cdots
\circ
\Theta _{B_{0} ', B_1 '}
$.

\end{lem}

{\bf Proof.}
Because
$(S_0 , B_0 )\blacksquare
\cdots \blacksquare
(S_m , B_m )=
(S_n ', B_n ')\blacksquare
\cdots \blacksquare
(S_0 ', B_0 ')
=(S_0 , B_0 )$
is a cycle,
and ${\cal S}$ is synchronized,
we have
$
\Theta _{B_{0} ', B_1 '} ^{-1}
\circ
\cdots
\circ
\Theta _{B_{n-1} ', B_n '}  ^{-1}
\circ
\Theta _{B_{m-1} , B_m }
\circ
\cdots
\circ
\Theta _{B_{0} , B_1 }
=
\Theta _{B_1 ' , B_{0} '}
\circ
\cdots
\circ
\Theta _{B_n ', B_{n-1} ' }
\circ
\Theta _{B_{m-1} , B_m }
\circ
\cdots
\circ
\Theta _{B_{0} , B_1 }
=
{\rm id} \hochdamit _{B_0 } ^{\Lambda  ^* (D_{S_0 } )\cdot S_0 }
$.
Hence
$
\Theta _{B_{m-1} , B_m }
\circ
\cdots
\circ
\Theta _{B_{0} , B_1 }
=
\Theta _{B_{n-1} ', B_n '}
\circ
\cdots
\circ
\Theta _{B_{0} ', B_1 '}
$.
\qed

\vspace{.1in}

By Lemma \ref{longthetawelldef}, the following is now well-defined.

\begin{define}

Let $(U,{\cal D} , G^* ;{\cal B} , {\cal F})$ be a
nested standard configuration.
Let the
saturated
subset ${\cal S}\subseteq {\cal F}$ be
synchronized.
For any two $(C,T), (S, B )\in {\cal S}$,
we define
$\Theta _{T,B} :=\Theta _{B_{m-1} , B_m }
\circ
\cdots
\circ
\Theta _{B_{0} , B_1 }
$,
where
$(C,T)=(S_0 , B_0 )\blacksquare
\cdots \blacksquare
(S_m , B_m )=
(S , B )$
is a path from $(C,T)$ to $(S,B)$.

\end{define}

Lemma \ref{mutualmaxlockcontrol} below now shows a
key property of
saturated synchronized sets
${\cal S}$:
Let $(S,B)\in {\cal S}$ and, as will happen in
a reverse depth enumeration,
let the action on $\Lambda  ^* (D)\cdot B$ be
fixed to be the identity.
Then, for all
$(C,T)\in {\cal S}$,
the action on
$B[\Lambda  ^* (D)\cdot S]$
together with
the image of $T$ determines the action on
$T[\Lambda  ^* (E)\cdot C]$.
Figure \ref{unsynch_factorial} shows that
${\cal S}$ being saturated is crucial:
Every $\blacksquare $-path from the bottom orbit
to the top orbit is synchronized, but the action on the bottom orbit
does not determine the action on the top orbit.

\begin{lem}
\label{mutualmaxlockcontrol}

Let $(U,{\cal D} , G^* ;{\cal B} , {\cal F})$ be a
nested standard configuration.
Let the
saturated
subset ${\cal S}\subseteq {\cal F}$ be
synchronized and let
$(C,T), (S, B )\in {\cal S}$.
Then, for all $\Phi , \Psi \in G^* $,
the following are equivalent.
\begin{enumerate}
\item
\label{mutualmaxlockcontrol1}
$\Phi [B]=\Psi [B]$
and
$
\Phi
\hochdamit
^{\Lambda  ^* (D_C )\cdot C }
_T
=
\Psi
\hochdamit
^{\Lambda  ^* (D_C )\cdot C }
_T $.
\item
\label{mutualmaxlockcontrol2}
$\Phi [T]=\Psi [T]$
and
$
\Phi
\hochdamit
^{\Lambda  ^* (D_S )\cdot S }
_B
=
\Psi
\hochdamit
^{\Lambda  ^* (D_S )\cdot S }
_B
$.
\end{enumerate}
In particular,
for all
$\Phi \in G^* _{\setcontrol S,T} $, we have
$\Phi
\hochdamit
^{\Lambda  ^* (D_C )\cdot C }
_T
=
{\rm id}
\hochdamit
^{\Lambda  ^* (D_C )\cdot C }
_T
$.

\end{lem}

{\bf Proof.}
Because
${\cal S}$ is saturated, either of
$\Phi [B]=B$
or
$\Phi [T]=T$
implies $\Phi [{\cal S} ]={\cal S}$.
Let
$(S_0 , B_0 )\blacksquare
(S_1 , B_1 )\blacksquare
\cdots \blacksquare
(S_m , B_m )$
be a path in ${\cal S}$.
Then
$\Theta _{B_0 , B_m } =
\Theta _{B_{m-1} , B_m }
\circ
\cdots
\circ
\Theta _{B_{0} , B_1 }
$.
Now, for any $\Phi \in G^* $ that maps ${\cal S}$ to itself,
we obtain the following by repeated
application of
Lemma \ref{thetatransmission}.
\begin{eqnarray*}
\lefteqn{
\Phi
\hochdamit
^{\Lambda  ^* (D_{S_m} )\cdot S_m }
\circ \Theta _{B_0 , B_m }
}\\
& = &
\Phi
\hochdamit
^{\Lambda  ^* (D_{S_m} )\cdot S_m }
\circ
\Theta _{B_{m-1} , B_m }
\circ
\cdots
\circ
\Theta _{B_{0} , B_1 }
\\
& = &
\Theta _{
\Phi \left[
B_{m-1} \right] , \Phi \left[
B_m \right] }
\circ
\cdots
\circ
\Theta _{\Phi \left[
B_{0} \right] , \Phi \left[
B_1 \right] }
\circ
\left.
\Phi
\hochdamit
^{\Lambda  ^* (D_{\Phi [S_0] } )\cdot \Phi [S_0 ] }
\right|
_{B_0 [\Lambda ^* (D_{B_0 } )\cdot S_0 ]}
\\
& = &
\Theta _{
\Phi \left[
B_{0} \right] ,
\Phi \left[
B_{m} \right]
}
\circ
\left.
\Phi
\hochdamit
^{\Lambda  ^* (D_{\Phi [S_0]} )\cdot \Phi [S_0 ]}
\right|
_{B_0 [\Lambda ^* (D_{B_0 } )\cdot S_0 ]}
,
\end{eqnarray*}

where
$\Theta _{
\Phi \left[
B_{0} \right] ,
\Phi \left[
B_{m} \right]
}
$
is well-defined, because $\Phi $
maps ${\cal S}$ to itself.

Now let $\Phi \in G^* $ such that
$\Phi [B]=B$ and
$
\Phi
\hochdamit
^{\Lambda  ^* (D_C )\cdot C }
_T
=
{\rm id}
\hochdamit
^{\Lambda  ^* (D_C )\cdot C }
_T
$.
Note that, in particular, this means $\Phi [T]=T$
and $\Phi [{\cal S} ]={\cal S}$.
Now, for every $\widetilde{S}\in
B [\Lambda  ^* (D_S )\cdot S ]$,
because
$\Theta
_{B,T}
\left(
\widetilde{S}
\right)
\in
T[\Lambda  ^* (D_C )\cdot C] $
we have
$
\Phi \left[
\Theta
_{B,T}
\left(
\widetilde{S}
\right)
\right]
=
{\rm id} \left[
\Theta
_{B,T}
\left(
\widetilde{S}
\right)
\right]
=
\Theta
_{B,T}
\left(
\widetilde{S}
\right)
$.
Hence, for every $\widetilde{S}\in
B [\Lambda  ^* (D_S )\cdot S ]$,
we have
$
\Phi \left[ \widetilde{S} \right]
=
\Phi \left[
\Theta _{T,B} \left(
\Theta _{T,B} ^{-1} \left(
\widetilde{S}
\right)
\right)
\right]
=
\Theta _{\Phi \left[ T\right] ,\Phi \left[ B\right] } \left(
\Phi \left[
\Theta
_{B,T}
\left(
\widetilde{S}
\right)
\right]
\right)
=
\Theta
_{T,B}
\left(
\Theta
_{B,T}
\left(
\widetilde{S}
\right)
\right)
=
\Theta _{T,B} \left(
\Theta _{T,B} ^{-1} \left(
\widetilde{S}
\right)
\right)
=
\widetilde{S}
$.
Hence
$
\Phi
\hochdamit
^{\Lambda  ^* (D_S )\cdot S }
_B
=
{\rm id}
\hochdamit
^{\Lambda  ^* (D_S )\cdot S }
_B
$.

This establishes
that \ref{mutualmaxlockcontrol1}
implies \ref{mutualmaxlockcontrol2}.
For the converse, exchange the roles of $(C,T)$ and $(S,B)$.

Finally, let
$\Phi \in G^* _{S,T} $.
Then $\Phi [T]=T$ and $\Phi
\hochdamit
^{\Lambda  ^* (D_S )\cdot S }
_B
=
{\rm id}
\hochdamit
^{\Lambda  ^* (D_S )\cdot S }
_B
$.
Consequently,
we have
$
\Phi
\hochdamit
^{\Lambda  ^* (D_C )\cdot C }
_T
=
{\rm id}
\hochdamit
^{\Lambda  ^* (D_C )\cdot C }
_T $.
\qed

\section{Outsized Factorial Components}
\label{maxlocksec}

Lemma \ref{disperseoutsized} below shows that
outsized factorial components are tailor-made for the type
of ``compensation by trading factors"
that
motivates Definition \ref{controlfunction}
and
that we have already seen in
Proposition \ref{breadth4Lemma}.
Unfortunately, Definition \ref{outsizedef} below formally
depends on the reverse depth enumeration
$\overrightarrow{{\cal B}} $.
We will see in
Proposition \ref{synchoutsizecond}
that this
is not an issue for synchronized components.
The fact that it is also not an issue for unsynchronized components
requires a substantial amount of detail and will be proved in \cite{SchAutCon2}.

\begin{define}
\label{lowerandupperlevels}

Let $(U,{\cal D} , G^* ;{\cal B} , {\cal F})$ be a
nested standard configuration.
For any subset ${\cal Q} \subseteq {\cal F}$,
the {\bf lower level} of ${\cal Q}$
is
$L({\cal Q} ):=\left\{ S:
(S,B)\in {\cal Q} \right\} $, and
the {\bf upper level} of ${\cal Q}$
is
$U({\cal Q} ):=\left\{ B:
(S,B)\in {\cal Q} \right\} $.
For any
$\blacksquare $-connected
${\cal C} \subseteq {\cal F}$
the {\bf vertex width} of ${\cal C}$
is the width of any factorial pair $(S,B)\in {\cal C}$.

\end{define}

\begin{define}
\label{outsizedef}

Let
$\left( U,{\cal D} , G^* ; \overrightarrow{\cal B}, {\cal F}\right)$
be
an opportunistically
enumerated nested
standard configuration.
Let
${\cal P}$ be a
saturated subset
of
$({\cal F},\blacksquare )$
of vertex width $w$.
We say
that
${\cal P}$ is
{\bf $\overrightarrow{\cal B} $-outsized} iff
there is an $n\in {\mat N}$
such that
$|U({\cal P})|\geq n\lg(w)$
and such that
$
\prod _{S\in {\cal B} \cap G^* \cdot L({\cal P}) }
\left|
\Lambda  ^* \left( \overleftarrow{\cal B} \not\setcontrol S \right)
\right|
\leq
(w!)^{n|\{ {\cal L}:{\cal L}\asymp {\cal P} \} |} $.
Otherwise, ${\cal P}$ is called
{\bf $\overrightarrow{\cal B}$-undersized}.

\end{define}

\begin{lem}
\label{disperseoutsized}

Let
$\left( U,{\cal D} , G^* ; \overrightarrow{\cal B}, {\cal F}\right)$
be
an opportunistically
enumerated nested
standard configuration
with control function $\beta _{\cal B} ^* $.
Let ${\cal P}$ be a
$\overrightarrow{\cal B} $-outsized
component of
$({\cal F},\blacksquare )$
such that,
for all $S\in {\cal B} \cap G^* \cdot L({\cal P})$,
we have
$\beta _{\cal B} ^* (S)
=
\left|
\Lambda  ^* \left( \overleftarrow{\cal B} \not\setcontrol S \right)
\right|
$.
Then
there is a control function
$\gamma _{\cal B} ^* $
such that
${\cal B} _u ^\gamma \subseteq {\cal B} _u ^\beta $ and
${\cal B} _u ^\gamma \cap G^* \cdot L({\cal P}) =\emptyset $.

\end{lem}

{\bf Proof.}
First note that we only need to
consider the case that there is an
$S\in {\cal B} \cap G^* \cdot L({\cal P})$
such that
$\beta _{\cal B} ^* (S)
>1$.

For
every
$S\in {\cal B} \cap G^* \cdot L({\cal P}) $,
we define
$\gamma _{\cal B} ^*  \left( S \right)
:=
2^{
{\lg \left( w!\right) \over w\lg (w)}
\cdot
|\Lambda  ^* (D_S )\cdot S |
}
$.
For all
other
$S\in {\cal B} $,
we define
$
\gamma _{\cal B} ^*  \left( S \right)
:=
\beta _{\cal B} ^*  \left( S \right)
$.
It is easily checked
by direct computation that, for $w\geq 5$, we have
${\lg \left( w!\right) \over w\lg (w)}
\leq d_w $
and
${\lg \left( w!\right) \over w\lg (w)}
\leq
0.99-{1.54 \over w }
$.
This establishes
part \ref{controlfunction2}
of Definition \ref{controlfunction}
for
$
\gamma _{\cal B} ^*
$
as well as,
via part \ref{level0elimlem1} of Proposition \ref{level0elimlem},
the claims about
${\cal B} _u ^\gamma $.
What remains to be checked is part \ref{controlfunction1}
of Definition \ref{controlfunction}.

Note that
$|L({\cal P} )|=w\cdot |U({\cal P}) |\geq w\lg(w)n$.
Now
\begin{eqnarray*}
\lefteqn{
\lg \left(
\prod _{S\in {\cal B} \cap G^* \cdot L({\cal P}) }
\left|
\Lambda  ^* \left( \overleftarrow{\cal B} \not\setcontrol S \right)
\right|
\right)
}\\
& \leq &
\lg \left( (w!)^{n|\{ {\cal L}:{\cal L}\asymp {\cal P} \} |} \right)
=
\lg \left( w!\right) n|\{ {\cal L}:{\cal L}\asymp {\cal P} \} |
=
{\lg \left( w!\right) \over w\lg (w)}
w
\lg \left( w\right) n
|\{ {\cal L}:{\cal L}\asymp {\cal P} \} |
\\
& \leq &
{\lg \left( w!\right) \over w\lg (w)}
|L({\cal P})| \cdot |\{ {\cal L}:{\cal L}\asymp {\cal P} \} |
=
{\lg \left( w!\right) \over w\lg (w)}
\left| G^* \cdot L({\cal P})| \right|
\\
& = &
{\lg \left( w!\right) \over w\lg (w)}
\sum _{S\in {\cal B} \cap G^* \cdot L({\cal P}) }
|\Lambda  ^* (D_S )\cdot S |
=
\sum _{S\in {\cal B} \cap G^* \cdot L({\cal P}) }
\lg \left( \gamma _{\cal B} ^*  \left( S \right)
\right)
\\
& = &
\lg \left(
\prod _{S\in {\cal B} \cap G^* \cdot L({\cal P}) }
\gamma _{\cal B} ^*  \left( S \right)
\right)
\end{eqnarray*}
and hence
\begin{eqnarray*}
\lefteqn{
\prod _{C\in {\cal B} }
\left|
\Lambda  ^* \left( \overleftarrow{\cal B} \not\setcontrol C \right)
\right|
}\\
& \leq &
\prod _{C\in {\cal B} }
\beta _{\cal B} ^* (C)
=
{
\prod _{C\in {\cal B} }
\beta _{\cal B} ^* (C)
\over
\prod _{S\in {\cal B} \cap G^* \cdot L({\cal P}) }
\beta _{\cal B} ^* (S)
}
\prod _{S\in {\cal B} \cap G^* \cdot L({\cal P}) }
\left|
\Lambda  ^* \left( \overleftarrow{\cal B} \not\setcontrol S \right)
\right|
\\
& \leq &
{
\prod _{C\in {\cal B} }
\gamma _{\cal B} ^* (C)
\over
\prod _{S\in {\cal B} \cap G^* \cdot L({\cal P}) }
\gamma _{\cal B} ^* (S)
}
\prod _{S\in {\cal B} \cap G^* \cdot L({\cal P}) }
\gamma _{\cal B} ^* (S)
=
\prod _{C\in {\cal B} }
\gamma _{\cal B} ^* (C)
\end{eqnarray*}

\vspace*{-.25in}
~\qed

\subsection{Synchronized Factorial Components}
\label{synchsubsec}

Lemma \ref{disperseoutsized} shows that the contribution of an
outsized component
of
$({\cal F},\blacksquare )$
(and its conjugates)
to the number of automorphisms can,
like the electron cloud in a Benzene ring,
be spread over all participating orbits rather than
kept localized in one orbit.
Unfortunately, Definition \ref{outsizedef}
formally depends on the
reverse depth enumeration $\overrightarrow{\cal B} $.
For synchronized components,
we can show
that being outsized
does not depend on
the particular
reverse depth enumeration $\overrightarrow{\cal B} $.

\begin{prop}
\label{synchoutsizecond}

Let $(U,{\cal D} , G^* ;{\cal B} , {\cal F})$ be a
nested standard configuration
and let ${\cal S}$ be a synchronized component of
$({\cal F},\blacksquare )$
of vertex width $w\geq 5$
such that $|U({\cal S})|\geq \lg(w)$.
Then, for every reverse depth enumeration
$\overrightarrow{\cal B}$,
we have that ${\cal S}$ is $\overrightarrow{\cal B}$-outsized.

\end{prop}

{\bf Proof.}
Let $\overrightarrow{\cal B}$ be
reverse depth enumeration
and let
${\cal B} \cap G^* \cdot L({\cal S})
=\left\{ S_{1 } , \ldots , S_{\ell } \right\} $,
enumerated such that $i<j$ implies that
$S_i $ occurs before $S_j $ in $\overrightarrow{\cal B} $.
Because ${\cal S}$ is synchronized,
by Lemma \ref{mutualmaxlockcontrol}
and Proposition \ref{folkloreonblocks2},
for all
$k\in \{ 2,\ldots , \ell \} $,
we have that
$
\left|
\Lambda  ^*
\left( \overleftarrow{\cal B} \not\setcontrol S_{k }  \right)
\right|
\leq
\left|
G^* _{\setcontrol S_{1} ,
S_{k } ^\dag }
\hochdamit ^{\Lambda  ^* (D_{k} )\cdot S_{k} }
\right|
=1$.
Moreover,
again because ${\cal S}$ is synchronized,
by Lemma \ref{mutualmaxlockcontrol},
we have that
$
\left|
\Lambda  ^*
\left( \overleftarrow{\cal B} \not\setcontrol S_{1 }  \right)
\right|
\leq
\left|
G^* _{\setcontrol S_{1 } ^\dag }
\hochdamit ^{\Lambda  ^* (D_{1} )\cdot S_{1} }
\right|
\leq
(w!)^{|\{ {\cal L}:{\cal L}\asymp {\cal S} \} |}$.
Now
\begin{eqnarray*}
\prod _{S\in {\cal B} \cap G^* \cdot L({\cal S}) }
\left|
\Lambda  ^* \left( \overleftarrow{\cal B} \not\setcontrol S \right)
\right|
& = &
\prod _{k=1} ^\ell
\left|
\Lambda  ^*
\left( \overleftarrow{\cal B} \not\setcontrol S_{k }  \right)
\right|
\leq
(w!)^{|\{ {\cal L}:{\cal L}\asymp {\cal S} \} |} ,
\end{eqnarray*}
which proves the claim.
\qed

\section{Main Results}
\label{mainres}

As indicated at the start of Section \ref{maxlocksec},
further work with unsychronized components
of $({\cal F} , \blacksquare )$
leads to a substantial
amount of new details that are best deferred to future work.
We now summarize the current state of affairs and provide a
brief outlook to this work.
Because of the large number of definitions that led us to this
point, references to the definitions used are given as superscripts.
First, we note that undersized components
of $({\cal F} , \blacksquare )$
are the only possible
cause for an automorphism set to have more than
$2^{0.99n} $ automorphisms.

\begin{theorem}
\label{upperboundonlyoutsizedleft}

Let
$\left( U,{\cal D} , G^* ; {\cal B} \right)$
be a nested\textsuperscript{\ref{primnestcoll}}
standard configuration\textsuperscript{\ref{standardconfdef}}
and let ${\cal U}$ be the set of
undersized\textsuperscript{\ref{outsizedef}}
components
of
$({\cal F},\blacksquare )$\textsuperscript{\ref{mutualmaxlockdef}}
of vertex width\textsuperscript{\ref{lowerandupperlevels}} $w\geq 5$.
Then
there is a
reverse depth enumeration\textsuperscript{\ref{RDEdef}}
$\overrightarrow{\cal B} $
of ${\cal B}$
with a control function\textsuperscript{\ref{controlfunction}}
$\gamma _{\cal B} ^* $
such that
$\displaystyle{
\left| G^* \right|
\leq
2^{0.99 |U|}
\prod _{D\in {\cal D} }
\Upsilon _{{\cal B} } (D)
} $, and
such that
${\cal B} _u ^\gamma
\subseteq L({\cal U}) $.

\end{theorem}

{\bf Proof.}
We start with an opportunistic reverse
depth enumeration
$\overrightarrow{\cal B}$
and a control function $\gamma _{\cal B} ^* $
from Theorem \ref{upperboundwithfactorials}.
Now, for every $\overrightarrow{\cal B} $-outsized
component of
$({\cal F},\blacksquare )$
of vertex width $w\geq 5$,
apply Lemma \ref{disperseoutsized}
to successively update
$\gamma _{\cal B} ^* $ until
${\cal B} _u ^\gamma
\subseteq L({\cal U}) $.
\qed

\vspace{.1in}

Therefore, clearly, if there are not too many points in
undersized components
of $({\cal F} , \blacksquare )$,
then the Automorphism Conjecture holds.

\begin{theorem}
\label{fewundersized}

The Automorphism Conjecture
holds for
the class ${\cal P}
$ of finite ordered sets
$P$
with $n$ elements
such that, for some choice of primitive nesting
collection\textsuperscript{\ref{primnestcoll}},
we have
$
|L({\cal U} )|\leq 0.005{n\over \lg (n)}
$.

\end{theorem}

{\bf Proof.}
First, let $P$ be an
ordered set
in this class
such that $P$ has $n$ elements and
does not contain
nontrivial order-autonomous antichains
and let $G^* :={\rm Aut} (P)$.
Pick a primitive nesting ${\cal B} $ such that
the set ${\cal U} $
satisfies the hypothesis.
By Theorem \ref{upperboundonlyoutsizedleft}
we have
\begin{eqnarray*}
|{\rm Aut} (P)|
& \leq &
2^{0.99 |U |}
\prod _{D\in {\cal D} }
\Upsilon _{{\cal B}  } (D)
\leq
2^{0.99 n}
\left| L({\cal U} )\right| !
\\
& \leq &
2^{0.99 n}
\left(
0.005{n\over \lg (n)}
\right) ^{
0.005{n\over \lg (n)}
}
\leq
2^{0.995 n} .
\end{eqnarray*}

Because, see \cite{Duffusprivate},
every ordered set has at least $2^n $ endomorphisms,
the Automorphism Conjecture holds for
the subclass of ${\cal P}$
in which no set contains nontrivial order-autonomous antichains.
By Proposition \ref{onelexiter},
the Automorphism Conjecture holds for
${\cal P}$.
\qed

\begin{remark}
\label{outlook}

{\rm
The author is quite certain that
all unsynchronized components are outsized,
and a proof will be included in
\cite{SchAutCon2}.
The natural next step for bounding the size of the automorphism group
is then to
compensate for undersized
(necessarily synchronized) components
${\cal C}$ of the factorial max-lock
graph.

Consider a synchronized undersized component
${\cal C}$
such that there are an $(S,B)\in {\cal C}$ and
a $(C,T)\in (G^* \cdot {\cal B} ^2 )\setminus {\cal C}$
such that $C\blacktriangleright ^B S$.
If $C$ occurs before $S$ and $S$
is the earliest block of $L({\cal C})$
in the chosen
reverse depth enumeration $\overrightarrow{\cal B}$,
then
$\left|
\Lambda  ^* \left( \overleftarrow{\cal B} \not\setcontrol S \right)
\right|
=1$,
and there is no contribution to $|G^* |$ from ${\cal C}$.

Otherwise, however, we must consider the case that
$C$ could control
blocks in the lower part of several
noncongruent
undersized components
${\cal C}_j $ of the factorial max-lock
graph.
This leads to structures similar to
what is used for unsynchronized components,
which will also be addressed in
the
follow-up \cite{SchAutCon2} to this work.
\qex
}

\end{remark}

\section{Proof of Concept: Width 9}
\label{takestock}

Theorems \ref{upperboundonlyoutsizedleft}
and \ref{fewundersized}
clearly show that, to
resolve the Automorphism Conjecture,
future work can focus on
the undersized
components of the
factorial max-lock
graph
$({\cal F},\blacksquare )$.
Indeed, if we can find a reverse depth enumeration
$\overrightarrow{\cal B}$ and a control function
$\gamma _{\cal B} ^* $
such that
${\cal B} _u ^\gamma =\emptyset $,
then
$\displaystyle{
|{\rm Aut}  (U)|
\leq 2^{0.99 |U|}
.
} $
Given the technical nature of these results,
the preceding sentence looks dangerously close to
the self-fulfilling
statement ``the conclusion holds when the theorem works."
To
demonstrate that
new conclusions can indeed be obtained
with the new tool without resorting to
self-fulfilling prophecies,
we will show how
Theorem \ref{upperboundonlyoutsizedleft}
can, with little effort, be used to establish the
Automorphism Conjecture for ordered sets of width at most $9$.
(Recall that the {\bf width} $w(P)$
of an ordered set $P$
is the size of the largest antichain contained in $P$.)
Ordered sets of bounded width
have recently gained attention in \cite{BonaMartin}.

The proof will exhibit a few more tools that
may be useful in an eventual resolution of the
Automorphism Conjecture.
We will forgo any and all attempts to
bargain with fate for an extension of
Theorem \ref{width9} below
from width $9$ to width $9+\Delta $,
and indeed we will present the argument for
Theorem \ref{width9} as rapidly as possible.
The natural next steps towards the
Automorphism Conjecture can only be to
successively
expand the
investigation of factorial pairs
until a resolution emerges.

\begin{define}

Let $(U,{\cal D} , G^* ;{\cal B},{\cal F} )$ be
a
nested standard configuration
and let
${\cal P}$ be a component of
$({\cal F},\blacksquare )$.
We call ${\cal P}$
{\bf thin} iff
$L({\cal P})$ consists of singletons.

\end{define}

\begin{lem}
\label{width9factcomp}

Let $(U,{\cal D} , G^* ;{\cal B},{\cal F} )$ be a
nested standard configuration
and let
${\cal P}$ be a component of
$({\cal F},\blacksquare )$ of vertex width $w\geq 5$.
If the width of $U$ is $9$ or less, then
${\cal P}$ is thin, synchronized, and
every $B\in U({\cal P})$ is an orbit.
Moreover,
if $B\in U({\cal P})$ is nontrivially woven with
an orbit $D$, then $D\in U({\cal P})$.

\end{lem}

{\bf Proof.}
If ${\cal P}$ were not thin, then the width of $U$ would
be at least $2w\geq 10$, which cannot be.
Hence ${\cal P} $ is thin.
Similarly, if $B\in U({\cal P})$
were not an orbit, then, because $B$ is a block, the containing
orbit would have at least
$2w\geq 10$ elements, which cannot be.
Hence
every $B\in U({\cal P})$ is an orbit.

Let
$(S,B)\in {\cal P}$
and let
$x\in U\setminus B$
such that $k:=|\{ s\in B: x\sim s\} |$
satisfies $k\in \{ 1, \ldots , |B|-1\} $.
In case
$k\in \{ 2, \ldots , |B|-2\} $,
we would have that
$|G^* \cdot x|\geq
|G^* \cdot \{ s\in B: x\sim s\} |\geq \pmatrix{|B|\cr k\cr } \geq 10 $,
which cannot be.
Thus $k\in \{ 1, |B|-1\} $.
Without loss of generality, we can assume that
$k=1$.
If any $b\in B$ was comparable to $2$ elements of $G^* \cdot x$, we
would obtain
$|G^* \cdot x |\geq 10$, which cannot be.
Therefore, for every $b\in B$, there is a unique
$x_b\in G^* \cdot x$
such that
$b\sim x_b $
and
$(S,B)$ and
$(\{ x\} ,G^* \cdot x)$
are max-locked.
Because $x\in U\setminus B$ was arbitrary, we conclude
in particular that,
if $B\in U({\cal P})$ is nontrivially woven with
an orbit $D$, then $D\in U({\cal P})$.

Because the width of $U$ is at most $9$, for any
two
$B^1 , B^2 \in U({\cal P})$
with $B^1 \not=B^2 $, there
are
$b^1 \in B^1 $ and $b^2 \in B^2 $
such that
$b^1 \sim b^2 $.
Therefore,
we can enumerate $U({\cal P})=:\left\{ B^0 , \ldots , B^h \right\} $
such that $i<j$ means there are
$b^i \in B^i $ and $b^j \in B^j $
such that
$b^i < b^j $.

For every $k\in \{ 0,\ldots , h\} $,
let
$B^k =:\left\{ x_1 ^k , \ldots , x_w ^k \right\} $.
Because ${\cal P}$ is $\blacksquare $-connected, for every
$k>0$ with $B^k \not= \emptyset $, we have
$B^{k-1} \not< B^k $.
Therefore, for every
$k>0$ with $B^k \not= \emptyset $,
and for every $j\in \{ 1, \ldots , w\} $,
we can assume that
$x_j ^{k-1} $ is either
the unique element of
$B^{k-1} $ that is below
$x_j ^k $
or $x_j ^{k-1} $ is
the unique element of
$B^{k-1} $ that is not below
$x_j ^k $.
Transitivity of the order now implies that,
for all $\ell <k$ we have $B^\ell  < B^k $, or,
for every $j\in \{ 1, \ldots , w\} $,
$x_j ^{\ell } $ is
the unique element of
$B^{\ell } $ that is below
$x_j ^k $
or it is
the unique element of
$B^{\ell } $ that is not below
$x_j ^k $.
In particular, we
obtain that ${\cal P}$ is synchronized.
\qed

\begin{theorem}
\label{width9}

The Automorphism Conjecture is true for
the class ${\cal W} _9 $ of ordered sets of width at most $9$.

\end{theorem}

{\bf Proof.}
Let $P$ be an ordered set of width at most $9$
with $n$ elements and no nontrivial order-autonomous
antichains.
Let
$\left( P,{\cal D} , {\rm Aut} (P) ; {\cal B} ,{\cal F} \right)$
be a nested
standard configuration.

By Lemma \ref{width9factcomp},
every component
${\cal P}$
of
$({\cal F},\blacksquare )$
is thin and synchronized.
Therefore, if ${\cal P}$
is undersized, then
$U({\cal P})\leq \lfloor \lg(9) \rfloor =3$.
Moreover, for every $w\in \{ 5,6,7,8,9\} $, we have that
$w!<2^{0.99\cdot 3w} $.
Therefore,
$|U({\cal P}|= 3$ implies
$\prod _{S\in {\cal B} \cap G^* \cdot L({\cal P}) }
\left|
\Lambda  ^* \left( \overleftarrow{\cal B} \not\setcontrol S \right)
\right|
\leq 2^{0.99|L({\cal P})|} $,
and hence, in this case,
we can
incorporate factors
$\left|
\Lambda  ^* \left( \overleftarrow{\cal B} \not\setcontrol S \right)
\right|
$
from $S\in L({\cal P})$ into the exponential
in
Theorem \ref{upperboundonlyoutsizedleft}.
Consequently,
because a component ${\cal P}$
of
$({\cal F},\blacksquare )$
with $|U({\cal P} )|=1$ would be a nontrivial order-autonomous antichain,
$
\prod _{D\in {\cal B} }
\Upsilon _{{\cal B} } (D)
$
is reduced to factors
$\left|
\Lambda  ^* \left( \overleftarrow{\cal B} \not\setcontrol S \right)
\right|
$
from components ${\cal Q}$
of
$({\cal F},\blacksquare )$
such that $|U({\cal Q})|=2$.
Let $f_P $ be the number of elements $x\in P$
that are contained in
the corresponding sets $L({\cal Q})$.

If
$f_P \leq \lg (n)$,
by Theorem \ref{upperboundonlyoutsizedleft},
we obtain
$|{\rm Aut} (P)|\leq
2^{0.99|P|} f_P !
\leq 2^{0.99n} \lg (n)!
\leq 2^{0.99n} \lg (n) ^{\lg (n)}
=2^{0.99n+ \lg (n) \lg( \lg (n)) }
$.

Otherwise, let
${\cal V}$ be the set of components
${\cal Q}$
of
$({\cal F},\blacksquare )$
that are thin, synchronized, and satisfy
$|U({\cal Q} )|= 2$,
and let
$V:=|{\cal V} |\geq {\lg(w)\over 18}.$
For each ${\cal Q}\in {\cal V}$,
let $w_{\cal Q}$ be the width of
the underlying set
$Q:=\bigcup L({\cal Q} )$, and note that
$Q$ is either the
disjoint union of $w_{\cal Q}$ chains with two elements each,
or
it is a set
$\left\{ x_1 ^0 , \ldots , x_{w_{\cal Q}} ^0 \right\}
\cup
\left\{ x_1 ^1 , \ldots , x_{w_{\cal Q}} ^1 \right\} $
with $x_j ^0 <x_i ^1 $ iff $j\not= i$ and no
further comparabilities.
In either case, it is easy to see that
there are at least $\left( w_{\cal Q} -1\right)^{w_{\cal Q}} $
endomorphisms of
$Q$ that map maximal elements to maximal elements and minimal elements to
minimal elements.

Similar to the proof of Proposition \ref{onelexiter},
the automorphisms of $P$ can be split into
the cosets of the normal subgroup $M$ that, for every
${\cal Q}\in {\cal V}$
maps each block in $U({\cal Q})$ to itself.
Every $[\Psi ]\in {\rm Aut} (P)/M$ contains exactly
$\prod _{Q\in {\cal V} } w_{\cal Q} !$
automorphisms.
By Lemma \ref{width9factcomp},
for every ${\cal Q} \in {\cal V}$,
no $B\in U({\cal Q})$ is nontrivially woven with
an orbit $D\not\in U({\cal Q})$.
Hence,
again similar to the proof of Proposition \ref{onelexiter},
for every $[\Psi ]\in {\rm Aut} (P)/M$,
there are at least
$\prod _{Q\in {\cal V} } \left( w_{\cal Q} -1\right)^{w_{\cal Q}} $
corresponding endomorphisms.
Thus
${|{\rm Aut}
(P)|
\over
|{\rm End}
(P)|}
\leq
{|{\rm Aut} (P)/M|
\prod _{Q\in {\cal V} } w_{\cal Q} !
\over
|{\rm Aut} (P)/M|
\prod _{Q\in {\cal V} } \left( w_{\cal Q} -1\right)^{w_{\cal Q}}
}
=
\prod _{Q\in {\cal V} }
{w_{\cal Q} !
\over
\left( w_{\cal Q} -1\right)^{w_{\cal Q}}
}
\leq
\left(
{5!\over 4^5 }
\right)
^{|{\cal V}|}
\leq
0.12^{\lg (n)\over 18}
$.

Therefore the
Automorphism Conjecture is true for
the class of ordered sets of width at most $9$
that do not contain nontrivial order-autonomous antichains.
The claim now follows from
Proposition \ref{onelexiter}.
\qed

\begin{remark}

{\rm
Example \ref{suspfactex}
and
Figure \ref{suspended_factorial}
show that, in general,
$({\cal F} , \blacksquare )$
can have components that are not thin.
In this case,
nontrivial weavings
with
primitive pairs that are not in ${\cal F}$
introduce
additional complications for any
attempts to
duplicate the arguments for thin components at the end of the
proof of Theorem \ref{width9}.
Aside from these possible complications
through non-controlling weavings with other
blocks,
Example \ref{noDendos} below shows another,
in the author's eyes more relevant,
complication, namely, that there are configurations
for which the argument
at the end of the
proof of Theorem \ref{width9}
is impossible.
\qex
}

\end{remark}

\begin{figure}

\centerline{
%\input{noproperendo.pic}
% This is a LaTeX picture output by TeXCAD.
% File name: [noproperendo.pic].
% Version of TeXCAD: 4.51
% Reference / build: 27-Nov-2018 (rev. a75)
% For new versions, check: http://texcad.sf.net/
% Options on the following lines.
%\grade{\on}
%\emlines{\off}
%\epic{\off}
%\beziermacro{\on}
%\reduce{\on}
%\snapping{\on}
%\pvinsert{% Your \input, \def, etc. here}
%\quality{8.000}
%\graddiff{0.005}
%\snapasp{1}
%\zoom{10.0000}
\unitlength 1mm % = 2.845pt
\linethickness{0.4pt}
\ifx\plotpoint\undefined\newsavebox{\plotpoint}\fi % GNUPLOT compatibility
\begin{picture}(102,17)(0,0)
\put(45,5){\circle*{2}}
\put(45,15){\circle*{2}}
\put(35,5){\circle*{2}}
\put(35,15){\circle*{2}}
\put(25,5){\circle*{2}}
\put(25,15){\circle*{2}}
\put(70,5){\circle*{2}}
\put(90,5){\circle*{2}}
\put(15,5){\circle*{2}}
\put(70,15){\circle*{2}}
\put(90,15){\circle*{2}}
\put(15,15){\circle*{2}}
\put(80,5){\circle*{2}}
\put(100,5){\circle*{2}}
\put(60,5){\circle*{2}}
\put(5,5){\circle*{2}}
\put(80,15){\circle*{2}}
\put(100,15){\circle*{2}}
\put(60,15){\circle*{2}}
\put(5,15){\circle*{2}}
\put(70,5){\line(0,1){10}}
\put(90,5){\line(0,1){10}}
\put(80,5){\line(0,1){10}}
\put(100,5){\line(0,1){10}}
\put(60,5){\line(0,1){10}}
\put(25,15){\line(-2,-1){20}}
\put(35,15){\line(-2,-1){20}}
\put(45,15){\line(-2,-1){20}}
\put(25,5){\line(-2,1){20}}
\put(35,5){\line(-2,1){20}}
\put(45,5){\line(-2,1){20}}
\put(5,5){\line(1,1){10}}
\put(15,5){\line(1,1){10}}
\put(25,5){\line(1,1){10}}
\put(35,5){\line(1,1){10}}
\put(15,5){\line(-1,1){10}}
\put(25,5){\line(-1,1){10}}
\put(35,5){\line(-1,1){10}}
\put(45,5){\line(-1,1){10}}
\put(45,15){\line(-4,-1){40}}
\put(45,5){\line(-4,1){40}}
\put(5,5){\line(3,1){30}}
\put(15,5){\line(3,1){30}}
\put(5,15){\line(3,-1){30}}
\put(15,15){\line(3,-1){30}}
\put(25,5){\oval(44,4)[]}
\put(25,15){\oval(44,4)[]}
\put(80,5){\oval(44,4)[]}
\put(80,15){\oval(44,4)[]}
\put(3,3){\makebox(0,0)[rt]{\footnotesize $D_1 $}}
\put(3,17){\makebox(0,0)[rb]{\footnotesize $D_2 $}}
\put(102,3){\makebox(0,0)[lt]{\footnotesize $D_3 $}}
\put(102,17){\makebox(0,0)[lb]{\footnotesize $D_4 $}}
%\emline(5,5)(60,15)
\multiput(5,5)(.1851851852,.0336700337){297}{\line(1,0){.1851851852}}
%\end
%\emline(5,15)(60,5)
\multiput(5,15)(.1851851852,-.0336700337){297}{\line(1,0){.1851851852}}
%\end
%\emline(15,5)(70,15)
\multiput(15,5)(.1851851852,.0336700337){297}{\line(1,0){.1851851852}}
%\end
%\emline(15,15)(70,5)
\multiput(15,15)(.1851851852,-.0336700337){297}{\line(1,0){.1851851852}}
%\end
%\emline(25,5)(80,15)
\multiput(25,5)(.1851851852,.0336700337){297}{\line(1,0){.1851851852}}
%\end
%\emline(25,15)(80,5)
\multiput(25,15)(.1851851852,-.0336700337){297}{\line(1,0){.1851851852}}
%\end
%\emline(35,5)(90,15)
\multiput(35,5)(.1851851852,.0336700337){297}{\line(1,0){.1851851852}}
%\end
%\emline(35,15)(90,5)
\multiput(35,15)(.1851851852,-.0336700337){297}{\line(1,0){.1851851852}}
%\end
%\emline(45,5)(100,15)
\multiput(45,5)(.1851851852,.0336700337){297}{\line(1,0){.1851851852}}
%\end
%\emline(45,15)(100,5)
\multiput(45,15)(.1851851852,-.0336700337){297}{\line(1,0){.1851851852}}
%\end
\end{picture}
}

\caption{A synchronized component of a factorial max-lock
graph such that all endomorphisms that
fix the orbits are automorphisms.}
\label{noproperendo}

\end{figure}
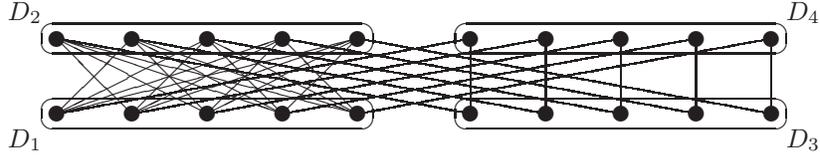

\begin{exam}
\label{noDendos}

{\rm
{\em
Let $M\in {\mat N}\setminus \{ 1,2\} $ and consider the
ordered set
$U=\bigcup {\cal D}$ with ${\cal D} =\{ D_1 , D_2 , D_3 , D_4 \} $,
$D_i =\{ x_1 ^i , \ldots , x_M ^i \} $ such that
$x_j ^1 <x_k ^2 $ iff
$j\not= k$,
$x_j ^3 <x_k ^2 $ iff
$j= k$,
$x_j ^3 <x_k ^4 $ iff
$j= k$,
$x_j ^1 <x_k ^4 $ iff
$j= k$, and no further comparabilities.
See Figure \ref{noproperendo} for a visualization.
Then every $f\in {\rm End}  (U)$
that maps every $D_i \in {\cal D}$ to itself
is an automorphism.
}

Let
$f\in {\rm End}  (U)$
map every $D_i \in {\cal D}$ to itself.
Because we are free to work with a power of $f$,
without loss of generality, we can assume that
$f^2 =f$
(see Proposition 4.4 in \cite{Schbook}).
Hence $f|_{f[U]} ={\rm id}|_{f[U]} $.

Let $x_j ^1 \in f[D_1 ]$.
Then
$f\left( x_j ^1 \right) =x_j ^1 $, which implies
$f\left( x_j ^4 \right) =x_j ^4 \in D_4 $, which implies
$f\left( x_j ^3 \right) =x_j ^3 \in D_3 $, which implies
$f\left( x_j ^2 \right) =x_j ^2 \in D_2 $.
Because $f\left( x_j ^2 \right) =x_j ^2 \in D_2 $,
no element of
$D_1 \setminus \{ x_j ^1 \} $ is mapped to
$x_j ^1 $.

Because $x_j ^1 \in f[D_1 ]$ was arbitrary,
no element in $D_1 \setminus f[D_1 ]$ can be mapped to any
element of $f[D_1 ]$.
Hence $f[D_1 ]=D_1 $ and, subsequently, for $j=2,3,4$,
$f[D_j ]=D_j $, which means
$f$ must be an automorphism.
\qex
}

\end{exam}

\vspace{.1in}

{\bf Acknowledgement.}
The author thanks Frank a Campo for a thorough
reading of an earlier, much more limited though
also rather technical, version of this paper,
for the detection of
mistakes, and many helpful suggestions;
and he thanks Nathan Bowler and his research group for
patiently listening to a presentation of an earlier stage of this work
and very useful discussions.

\end{document}